\documentclass[12pt]{article}
\usepackage{amsfonts,amssymb,version}
\def\A{\mbox{$\mathbb A$}}
\def\Q{\mbox{$\mathbb Q$}}
\def\P{\mbox{$\mathbb P$}}
\def\C{\mbox{$\mathbb C$}}
\def\F{\mbox{$\mathbb F$}}
\def\G{\mbox{$\mathbb G$}}
\def\E{{\cal E}}
\def\Cu{{\cal C}}
\def\sA{\mbox{\scriptsize{$\A$}}}

\def\mf{\mathfrak}

\def\nin{{\not\in}}

\def\Z{\mbox{$\mathbb Z$}}

\let\hra\hookrightarrow

\let\ra\rightarrow

\let\lra\longrightarrow

\let\ov\overline

\let \wt \widetilde

\def\vequal{\mbox{\large $\parallel$}}
\newcommand{\liminv}[1]{{\displaystyle{\mathop{\rm lim}_{\buildrel\longleftarrow\over{#1}}}}\,}
\newtheorem{theorem}{Theorem}[section]
\newtheorem{proposition}[theorem]{Proposition}
\newtheorem{lemma}[theorem]{Lemma}
\newtheorem{definition}[theorem]{Definition}
\newtheorem{corollary}[theorem]{Corollary}

\def\rem{\refstepcounter{theorem}\paragraph{Remark \thetheorem}}

\def\proof{\paragraph{Proof}}

\makeatletter \def\l@section{\@dottedtocline{1}{0em}{1.2em}} \makeatother

\begin{document}
\baselineskip=15pt

\title{Parabolic reductions of principal bundles}

\author{Yogish I, Holla}
\date{April 17, 2002}
\maketitle

\begin{abstract} 
In this paper we describe the structure of the space of parabolic reductions, 
and their compactifications, of principal $G$-bundles over a 
smooth projective curve over an algebraically closed field of 
arbitrary characteristic. We first prove estimates for the 
dimensions of moduli spaces of stable maps to the twisted flag varieties 
$E/P$ and Hilbert schemes of closed 
subschemes of $E/P$ with same Hilbert polynomials as that of a $P$-reductions
of $E$. This generalizes the earlier results of Mihnea Popa and 
Mike Roth to connected reductive groups and the results of 
Y. I. Holla and M. S. Narasimhan to the case of non-minimal sections.
We then prove irreducibility and generic smoothness of the space of 
reductions for large numerical constraints, using the above result and 
the methods of G. Harder. 
We also study these space in more detail for generically stable $G$-bundles. 
As a consequence we can generalize the lower bound results of H. Lange to $G$.   
\end{abstract}
\section{Introduction}
Let $C$ be a smooth projective over an algebraically closed field $k$ 
of arbitrary characteristic. Let $G$ be a connected reductive algebraic 
group. Let $P$ be a parabolic subgroup. Let $E$ be a principal $G$-bundle over $C$. The main objective of this paper is 
to understand the structure of the
space of $P$-reductions of $E$ and its compactifications. 
Let $\pi: E/P \ra C$ be the associated $G/P$ bundle over $C$. 
Let $T_{\pi}$ be the tangent bundle along the fibers of the morphism $\pi$. 
For a $P$ reduction  $\sigma$ of $E$, equivalently a section of $\pi$,
we will denote by $T_{\sigma} =\sigma ^*(T_{\pi})$ the normal bundle of 
$\sigma(C)$. Recall that a section $\sigma$ is minimal if 
${\rm deg}(T_{\sigma})$ is minimal among the degrees of the normal bundles 
for all sections of $\pi$. 
In Holla-Narasimhan \cite{Holla-Narasimhan} it was proved that 
the every irreducible component of the Hilbert scheme of closed subschemes 
of $E/P$ containing a minimal section has a dimension bound of 
${\rm dim}(G/P)$.

For a $P$ reduction $\sigma$ we define its numerical type 
$[\sigma] \in {\cal X}_*(P)={\rm Hom}({\cal X}^*(P),\,\Z)$ by 
$[\sigma](\chi)={\rm deg}(E_{\sigma}\times ^P \chi )$, where 
$E_{\sigma}\times ^P \chi $ is the line bundle associated to the $P$-bundle 
$E_{\sigma}$ defined by $\sigma$ via the character $\chi $ of $P$.

Fixing a polarization of the curve $C$ we can define a Hilbert scheme
${\rm Hilb}^{[\sigma]}_{E/P}$ which parameterizes closed subschemes of $E/P$ 
whose Hilbert polynomials with respect to a generating set of polarizations of $G/P$ coincide with that of $[\sigma]$. We also have the open subscheme 
${\rm Sec}^{[\sigma]}_{E/P}$ of the above Hilbert scheme parameterizing the 
space of sections of $\pi$. We have a partial ordering on the set 
${\cal X}_*(P)$ defined by 
$[\sigma _1] \le [\sigma _2]$ if for every dominant character $w$ of
$P$ we have $([\sigma _2] - [\sigma _1],\,w) \in \Z^{(\ge 0)}$ and
$([\sigma _2] - [\sigma _1],\,\chi)=0$ for every $\chi \in {\cal
X}(G)$.  This defines a notion of a numerically minimal sections (types) 
that is those section (types) for which $[\sigma]$ is maximal with respect the
above ordering among the numerical types of $P$ reductions of $E$.
We show that given a $G$-bundle there are only finitely many minimal 
numerical types and for all these types we have similar dimension bounds 
as in the case of \cite{Holla-Narasimhan}. 
Our first result is a generalization of the dimension estimates for the Hilbert schemes  in the case of non minimal sections.  

Let $\gamma _1,\, \ldots,\, \gamma _m$ be 
the set of minimal numerical types for $E$. We show that 
if $X$ is an irreducible component of ${\rm Hilb}^{[\sigma]}_{E/P}$ which 
contains the reduction of structure group $\sigma$ as a Hilbert point
then there exists an $i\in \{1,\,\ldots,\,m\}$ such that 
$[\sigma] \le \gamma _i$ and 
$${\rm dim}(X) \le {\rm dim}(G/P)+
d({[\sigma]})-d({\gamma_i}).$$
The above result is proved using a similar dimension estimates for the 
moduli space of maps. Let ${\ov M}_g(E/P,\beta_{[\sigma]})$ be the moduli 
space of stable maps from genus $g$ curves to $E/P$ with 
$\beta _{[\sigma]} \in H^2(E/P,\,\Q_l)^*$ be a class determined by the 
numerical type and the fixed polarization of $C$. 
 We show that 
if $X$ is an irreducible component of ${\ov M}_{g}(E/P, \beta_{[\sigma]})$
then there exists an $i\in \{1,\,\ldots,\,m\}$ such that 
$[\sigma]\le \gamma _i$
$${\rm dim}(X) \le {\rm dim}(G/P)+
d([\sigma])-d(\gamma_i).$$
The above result is a generalization of a similar result of 
Mihnea Popa and Mike Roth \cite{Popa-Roth} for the case $G=GL_n$ and $P$ a 
maximal parabolic and the method of proof is similar to this case. 

We next address the question of generic smoothness and irreducibility of 
the space of sections of $\pi$. We fix a root system of $G$ by considering a
Borel subgroup $B$ and a maximal torus $T$. 
We say a numerical type 
$[\sigma] \in{\cal X}_*(P)$ satisfies the property ($*$) for $N$ if 
$[\sigma](\chi) \le -N$ for every non-trivial character of $P$ which when 
restricted to the maximal torus $T$ is a non-negative linear combination of 
simple roots with respect to the root system. 
We show that there exists an integer $N$ such that 
if $E$ admits a $P$ reduction of numerical type $[\sigma]$ satisfying the 
property $(*)$ for $N$ then 
${\rm Sec}^{[\sigma]}_{E/P}$ is irreducible and generically smooth 
of expected dimension $d([\sigma])+(1-g){\rm dim}(G/P)$.  
This result is proved for the case $G=GL_n$ and $P$ a maximal parabolic 
in \cite{Popa-Roth} and for arbitrary $G$ and $P=B$ a Borel subgroup when 
the curve $C$ is over a finite field in Harder \cite{Harder}. Using the 
methods in \cite{Harder} we first derive the above result for Borel subgroups 
over arbitrary fields and then generalize this to case of parabolic subgroups.

Next we define the notion of a generically stable $G$-bundles extending a 
similar notion for vector bundles as defined in  Example 5.7 of 
\cite{Popa-Roth}. Our main result is the existence of generically 
stable bundles when the genus of the curve $C$ is atleast two. 
We also study some basic properties 
enjoyed by these $G$-bundles. For example we show that if $E$ is a 
principal $G$ bundle which is generically stable then 
$d([\sigma])\ge (g-1){\rm dim}(G/P)$
for every $P$-reduction $\sigma$. This generalizes the lower bound result of 
Lange (see \cite{Lange2}). Some of these ideas should also lead to 
the computation of Gromov-Witten invariants for the twisted flag varieties $E/P$ and this will be done elsewhere. 

For a fixed $G$-bundle $E$ we have two different compactifications namely the 
Hilbert scheme and the moduli space of maps. In general one can show that there
are no morphisms between them. Also there are several partial Drinfeld 
compactifications we can define using representations of $G$. It is possible 
to understand which of these are the images of morphisms from earlier 
compactifications (see Gaitsgory-Braverman \cite{Gaitsgory-Braverman} for 
an account on Drinfeld compactifications). This part will appear elsewhere.

In proving these results we need several basic facts about principle bundles 
over curves. In Section 2, we prove a technical result about existence of 
$B$-reductions of 
$E$ following ideas of Ramanathan \cite{Ramanathan83} and as a 
consequence we recover several basic properties of principal $G$-bundles.  
We partially 
answers a question of Friedman-Morgan on the behavior strata defined by 
the Harder-Narasimhan reduction. This is in the case when such a reduction 
is defined by Borel subgroups. 

For a principal $G$-bundle $E$ we define a 
canonical element $c(E)\in {\cal X}_*(T)/{\hat Q}$, where ${\hat Q}$ is the 
coroot lattice of $G$, using a Borel reduction and show that this element 
$c(E)$ exactly parameterizes the algebraic equivalence classes of $G$-bundles 
over $C$, thus generalizing the fact that fundamental group of $G$ 
parameterizes the topological equivalence classes of $G$-bundles over 
$\C$ to arbitrary characteristic. 

For an element $c \in {\cal X}_*(T)/{\hat Q}$ and a positive integer $d$,
let $M(c,d)$ be the set of isomorphism classes of $G$-bundles  
of topological type $c$ such that the instability degree 
${\rm Ideg}(E) = {\rm Max} \{ {\rm deg}({\rm ad}E_{\sigma}) |
(P,\sigma)\} \le d$, where maximum is taken over all $P$ reductions $\sigma$
of $E$ and all parabolic subgroups. We construct a finite type
 irreducible smooth scheme $S$ and a family of $G$-bundles $\E$ over 
$C\times S$ which is versal at every $x\in S$ and such that for 
each $x\in S$ the bundle $\E_{C\times \{x\}}$ lies 
in $M(c,d)$ and every member in $M(c,d)$ is an occurs in the above family.
We also show existence of stable bundles for curves of genus 
atleast two.

This paper is organized as follows. In Section 2 we describe the numerical 
types of $P$-reductions and study some basic properties. Section 3 deals 
with  Borel reductions and here we prove results on algebraic equivalence 
and versal families and existence of stable bundles. The dimension 
estimates for the Hilbert schemes and moduli space of stable maps is dealt in
Section 4. The irreducibility and generic smoothness is the content of Section 
5. Finally in the last Section we prove results about generic stability of 
$G$-bundles.

{\bf Acknowledgments} I would especially like to thank M.S. Narasimhan 
for valuable suggestions and encouragement.
I would also like to thank G. Harder for, among other things, 
explaining the proof of his Theorem and how it might be extended to 
arbitrary fields. 
I thank  B. Fantechi,
L. G\"ottsche, N. Nitsure  and  Mike Roth for some useful 
conversations. I would like to acknowledge  
Max-Planck Institute, Bonn and 
Abdus Salam International center for Theoretical Physics, Trieste for 
hospitality where the work was done.

\section{Some basic facts about principal bundles}
In this section we recall and prove some basic facts about principal $G$ 
bundles on $C$. 
Let $k$ be an algebraically closed field. 
Let $G$ be a connected
reductive algebraic group over $k$. Let $T$ be a maximal torus and $B$ a Borel
subgroup containing $T$. Let $U$ be the unipotent radical of $B$.
Then $B$ is a semi-direct product $U\cdot T$. Let $i :T\hra B$ and $j:B
\hra G$ be the inclusions. and $p_B: B\lra B/U=T$ be the projection.
Let $W=N(T)/T$ be the Weyl group and $w_0 \in W$ the element of
maximal length in $W$.  Let $\G _m$ be the one dimensional torus and
$\G _a$ the additive group.  We denote by ${\cal X}_*(T)$ be the group
of 1-parameter subgroups of $T$ (denote by 1-PS). ${\cal X}^*(T)$
denotes the group of characters of $T$.  We have a perfect pairing
${\cal X}_*(T) \otimes {\cal X}^*(T) \lra \Z$ which will be denoted by
$(\cdot,\,\cdot)$.  Let $\Phi \subset {\cal X}^*(T)$ be the root system of
$G$, $\Phi ^+$ be the set of positive roots and $\Delta\,=\
\{ \alpha_1,\,\ldots,\,\alpha _r\}$ the set of simple roots
corresponding to $B$. For any $\alpha\in \Phi$, let $T_{\alpha}$ the 
connected component
of ${\rm ker}(\alpha)$ and $Z_{\alpha}$ the centralizer of
$T_{\alpha}$ in $G$. Then the derived group
$[Z_{\alpha},\,Z_{\alpha}]$ is of rank one and there is a unique 1-PS
${\hat \alpha}: \G _m \ra T\cap [Z_{\alpha},\,Z_{\alpha}]$ such that
$T=({\rm im}\,{\hat \alpha}) \cdot T_{\alpha}$ and 
$({\hat \alpha},\,\alpha)=2$. This ${\hat \alpha}$ is the coroot corresponding to
$\alpha$. We denote by ${\hat \Phi}$ the set of coroots. The quadruple
$\{{\cal X}^*(T),\, \Phi,\,{\cal X}_*(T),\,{\hat \Phi}\}$ defines a
root system. For each $\alpha \in \Phi$ we have the fundamental dominant 
weight $w_{\alpha}\in {\cal X}^*(T)\otimes \Q$ defined by 
$({hat \beta},w_{\alpha})=\delta _{\alpha, \beta}$ and 
$(\gamma,w_{\alpha})=0$ for any 1-parameter group in the connected component of the center of $G$.
Let $Q \subset {\cal X}^*(T)$
 (resp.${\hat Q} \subset{\cal X}_*(T)$)  be the (co)-root lattice generated 
by $\Phi$ (resp ${\hat \Phi}$).
We have a partial ordering $\le$ in ${\cal X}_*(T)$ defined by 
$\mu \le \lambda$ if and only if $(\lambda -\mu,\,w_{\alpha}) \in \Z ^{(\ge 0)}$ and 
$(\lambda -\mu,\,\chi)=0$ for $\chi \in {\cal X}(G)$.

Let $P$ be a parabolic subgroup of $G$ containing $B$. Let $R_u(P)$ be its 
unipotent radical.
then there is a subset $I \subset \Delta$ such that $P=P_I$
Let $Z_I=(\cap _{\alpha \in I}{\rm ker} \alpha)^0$ be the connected
component of the intersection of the kernels of the roots in $I$. then
we have a chosen Levi decomposition $P=R_uP \cdot L$ such that 
$L$ a Levi subgroup containing
$T$ defined to be the centralizer of $Z_I$.
We will fix such a splitting $i : L\lra P$. 

Let $C$ be an smooth projective curve over $k$.  
Let $E$ is a principal $G$ bundle over $C$.  
Let $\sigma$ be a 
reduction of structure group of $E$ to $P$. By this we mean a pair
$\sigma = (E_{\sigma},\phi)$ with $E_P$ a principal $P$-bundle and an 
isomorphism $\phi:E_P \lra E$, equivalently a reduction of structure group is 
 a section $\sigma$ of the fiber bundle $\pi: E/P \lra C$.  
Here $E/P$ denotes the extended fiber bundle $E \times ^G G/P$ over $C$.

Let $T_{\pi}$ be the tangent bundle along the fibers of the map $\pi$. 
For a reduction of structure group $\sigma$ we will denote by $T_{\sigma}$ 
the vector bundle defined by the pull back of $T_{\pi}$ under $\sigma$. 
We will also fix notations for the Lie algebras by putting ${\mathfrak
g}$, ${\mathfrak p}$, ${\mathfrak m}$, ${\mathfrak u}$ for Lie
algebras of $G$, $P$, $L$ and $R_uP$ respectively.
Then we see that $T_{\sigma}$ is the vector bundle on $C$ 
associated to $E_{\sigma}$ for the representation of $P$ on ${\mf g}/{\mf p}$.

First we state a lemma which bounds the degree of the tangent bundle along the fibers of the map $\pi$. 

\begin{lemma}\label{minimal}
There exists a constant $C$ (independent of $\sigma$) 
such that for any reduction $\sigma$ of $E $ to $P$ we have 
${\rm deg}(T_{\sigma})\ge C$
\end{lemma}
\proof (see Lemma 2.1, \cite{Holla-Narasimhan}).$\hfill \square$

Recall that the above lemma enables us to define the notion of a minimal 
reduction namely those reductions of structure groups for which 
${\rm deg}(T_{\sigma})$ is minimal.

There is a stronger notion of minimality of reductions we will be 
interested. Suppose $E_{\sigma}$ is the $P$-bundle
associated to $\sigma$ then we define an element $[\sigma]$ of
${\cal X}_*(P)={\cal X}_*(L) = {\rm Hom}({\cal X}^*(L),\, \Z)$ by assigning 
$[\sigma](\chi)= {\rm deg}(\chi _*(E_{\sigma}))$.  This point of
${\cal X}_*(L)$ actually determines the Atiyah-Bott point in the sense
of Friedman-Morgan \cite{Friedman-Morgan1} of the reduction $\sigma$.
We say that $[\sigma]$ is the {\bf numerical type} of the reduction $\sigma$.

Now we can define a partial ordering on the elements of
${\cal X}_*(L)$ as follows. We say the numerical types 
$[\sigma _1] \le [\sigma _2]$ if for every dominant character $w$ of
$P$ we have $([\sigma _2] - [\sigma _1],\,w) \in \Z^{(\ge 0)}$ and
$([\sigma _2] - [\sigma _1],\,\chi)=0$ for every $\chi \in {\cal
X}(G)$.  We say a reduction of structure group $\sigma$ is 
{\bf numerically minimal}
if the numerical type $[\sigma]$ is maximal with respect the
above ordering among the numerical types of $P$ reductions 
$E$.  Here one observes that the second
condition in the definition of the partial ordering is automatically
satisfied if the the numerical types corresponds to P reduction 
of a fixed principal $G$-bundle $E$.

One notes that ${\rm deg}(T_{\sigma})$ depends only on the numerical type
$[\sigma]$ corresponding to a reduction $\sigma$ and not on the
reduction itself hence we denote $d([\sigma])={\rm deg}(T_{\sigma})$ for
some section $\sigma$ whose numerical type is $[\sigma]$. 
Also one observes the definition of $d([\sigma])$ can be enlarged to 
define it for any element of ${\cal X}_*(L)$ by 
$d([\sigma])=[\sigma](\chi _P)$, where $\chi _P$ is the character of $P$,
hence $L$, defined by the highest exterior power of the representation 
of $P$ on ${\mf g}/{\mf p}$.

\rem A reduction of structure group  $\sigma$ is numerically minimal if 
it is minimal with respect to the degrees of the tangent bundle $T_{\sigma}$. 
More generally $\sigma$ is 
minimal if ${\rm deg}(\sigma ^*(L))$ is minimal with respect to all reductions
for a fixed line bundle $L$ which is ample along the fibers of $\pi$ 
(this corresponds to negative powers of dominant character of $P$ upto 
tensoring with a line bundle pulled up from $C$).

The first basic lemma is the existence of numerically minimal sections. 

\begin{lemma}\label{lowerbound}For a given $N$,
there are only finitely many numerical types 
$[\tau]$, defined by $P$-reductions of a fixed $G$-bundle $E$, 
with the property that 
$d([\tau])\le N$. Moreover given any $P$-reduction $\sigma$
there exists an numerically minimal reduction $\sigma _0$ such that  
$[\sigma] \le [\sigma _0]$.
\end{lemma}
\proof We first prove the lemma for the case when the parabolic $P$ is maximal.
Further we may assume that the parabolic $P$ is of the form 
$P_{\Delta -\alpha}$
for some $\alpha \in \Delta$ by fixing all the root datum.
Now both parts of the lemma \ref{lowerbound} follows from the 
Lemma \ref{minimal} and the fact that $[\sigma]$ is completely determined by
$d([\sigma])$. This proves the lemma in this case.
For a given $\alpha \in \Delta$ let $n_{\alpha}$ be a positive integer such 
that $-n_{\alpha}w_{\alpha}$ defines the character 
$\chi _{P_{\Delta -\alpha}}$.
For the general case, we assume the parabolic $P$ is of the form 
$P=P_I$ for a subset $I\in \Delta$.  
The numerical type of $\tau$ is completely determined by what
values it takes on the characters $n_{\alpha}w_{\alpha}$ for 
$\alpha \in I$ and the characters of the
group $G$.   
Since any reduction to $P$ automatically determines a reduction to the 
parabolic $P_{\Delta -\alpha}$, and using the conclusion of the lemma for 
the maximal
parabolics, we see that $[\sigma](n_{\alpha} w_{\alpha})$  is
bounded above as we vary $\sigma$ over the reductions to $P$.
Hence we see that subset of ${\cal X}_*(L)$ we are interested in 
is finite.
The second assertion in the lemma is a consequence 
of the first. $\hfill \square$

Recall that a principal $G$-bundle is said to be (semi)stable if for any 
reduction of structure group $\sigma$ of $E$ to any parabolic we have 
${\rm deg}(T_{\sigma}) >(\ge)0$. This definition is equivalent to the condition
that for any maximal parabolic $P$ and a dominant character $w$ of $P$ we have 
${\rm deg}(w_*(E_{\sigma}))< (\le) 0$. 

If $P$ is a Borel subgroup $B$ then the numerical type $[\sigma]$ defines a 
1-PS on the maximal torus $T$.
Further for any $T$-bundle $E_T$ we can similarly define its numerical type 
by a 1-PS $[E_T]$ defined by $[E_T](\chi)={\rm deg}(\chi _*(E_T))$.

Fix a polarization of the curve $C$.
Recall that for a principal $G$ bundle ${\cal E}$ over $C\times S$ for 
a scheme $S$ of finite type over $k$ there is a projective scheme 
${\rm Hilb}^p_{{\cal E}/P,S}$ over $S$ parameterizing 
the closed subschemes of ${\cal E}/P$ flat over $S$ 
with a fixed Hilbert polynomial $p$.
When the parabolic is not maximal we see that the Hilbert scheme
further decomposes into open and closed subschemes owing to the fact
that we have many polarizations of $G/P$.
Hence we have to fix a finite set of generating  polarizations of $G/P$
 to set a topological type. One can check that these polarizations can be computed by the numerical type $[\sigma]$ of reductions.
Hence we define the Hilbert scheme ${\rm Hilb}^{[\sigma]}_{{\cal E}/P}$ 
which 
parameterizes closed subschemes of  ${\cal E}/P$ flat over $S$ with all
Hilbert polynomials same as that of $\sigma$. This defines an open and closed 
subschemes of the Hilbert scheme defined above.
Also we have an open subscheme of the Hilbert scheme 
${\rm Hilb}^{[\sigma]}_{{\cal E}/P}$ corresponding to the subschemes which are actually sections of the morphism $\pi_S:\E/P \ra S$. 
We will denote this scheme by
${\rm Sec}^{[\sigma]}_{{\cal E}/P}$.
One of the properties 
 of the Hilbert scheme and the space of sections that will be useful to us 
is their behavior under the base change. Namely for any morphism 
$S'\lra S$ we have natural isomorphisms
$ {\rm Hilb}^{[\sigma]}_{{\cal E}_{S'}/P}\cong 
{\rm Hilb}^{[\sigma]}_{{\cal E}/P}\times _SS'$ and 
${\rm Sec}^{[\sigma]}_{{\cal E}_{S'}/P}\cong
{\rm Sec}^{[\sigma]}_{{\cal E}/P}\times _SS'$. Here ${\cal E}_{S'}$ denotes the pull back of ${\cal E}$ under the morphism $S'\ra S$.

we now prove a Lemma about the space of sections which will be used later.
\begin{lemma} \label{lclosed}
Let $[\sigma]$ be an numerical type. Let $\E$ be a family of $G$-bundles 
over $C\times S$ with $S$ a finite type scheme. Then the subset of points in 
$S$ corresponding to principal $G$-bundles which admits $P$ reduction of
numerical type $[\sigma]$ is constructible. 
\end{lemma}
\proof 
The lemma follows from the fact that the subset we are interested in is exactly the image of the morphism 
${\rm Sec}^{[\sigma]}_{{\cal E}/P}\ra S$. 
Since the space of sections are 
finite type over $S$ hence the image is constructible.
$\hfill \square$

Let $P_1 \subset P$ be two parabolic subgroups of $G$. Let $L$ (resp. $L_1$ 
be the Levi quotients of $P$ (resp. $P_1$) and 
$Z_0(L)$ (resp. $Z_0(L_1)$) be the connected component of the center.
Let ${\ov L}$ denote the quotient $L/Z_0(L)$. Let ${\ov P_1}$ be the parabolic
subgroup of ${\ov L}$ defined by the image of   
$P_1$ under the natural map $p:P\ra {\ov L}$.
Let $E$ be a principal $G$-bundle over $C$ and 
let $\sigma$ be a $P$ reduction of 
$E$. Let $E_1=p_*(E_{\sigma})$ be the associated 
principle ${\ov L}$-bundle over $C$.
We have the following lemma.
\begin{lemma}\label{compare} 
If $[{\ov \sigma}_1] \in {\cal X}_*({\ov P}_1)$ is such that 
$E_1$ admits a  ${\ov P}_1$ reductions of numerical type 
$[{\ov \sigma}_1]$ then there is a unique 
$[\sigma _1] \in {\cal X}_*(P_1)$ such that there is a bijective correspondence between the $P_1$ reductions of 
$E_{\sigma}$ of numerical type $[\sigma _1]$ and 
${\ov P}_1$ reductions of $E_1$ of numerical type $[{\ov \sigma}_1]$. 
\end{lemma} 
\proof 
The main observation is that we have natural isomorphisms 
$P/P_1\cong {\ov L}/{\ov P}_1$. From here it follows that there is a natural 
bijection between  $P_1$ reductions of $E_{\sigma}$ and 
${\ov P}_1$ reductions of $E_1$. 
Let ${\ov \sigma}_1$ be a reduction of 
structure group of $E_1$ to ${\ov P}_1$ and let $\sigma_1$ be the $P_1$ 
reduction of $E_{\sigma}$ which corresponds to ${\ov \sigma}_1$ under the above bijection. Now the natural morphism
 ${\cal X}_*(P_1)\ra {\cal X}_*(P_1)$ takes $[\sigma _1]$  to $[\sigma]$ and 
the natural morphism ${\cal X}_*(P_1)\ra {\cal X}_*({\ov P}_1)$ 
takes $[\sigma _1]$  to $[{\ov \sigma}_1]$. The proof of the lemma will be complete once we establish the injectivity of the homomorphism
${\cal X}_*(P_1)\ra {\cal X}_*(P_1)\oplus {\cal X}_*({\ov P}_1)$. The last statement follows from the fact that any character of $P_1$ can be uniquely written 
as a rational linear combination of a character of $P$ and a character of 
${\ov P_1}$. $\hfill \square$

\section{Borel reductions and algebraic equivalence}

In this section prove a result about Borel reductions of $G$-bundles and 
as a consequence we derive results on algebraic 
equivalence of principal bundles and irreducibility of the moduli spaces etc.

We will now define a notation which will be used through in the article.
We say a point $[\sigma] \in \chi _*(L)$(or a reduction of structure group 
$\sigma$) satisfies the  
{\bf property ($*$) for} $N$  if 
\begin{eqnarray*}
-[\sigma](\chi ) & \ge N & \mbox{for every non-trivial character}~\chi~
\mbox{which 
when restricted to}~T \\
& &  \mbox{is a non-negative linear combination of simple roots}. 
\end{eqnarray*}
Note that a reduction $\sigma$ of $E$ to the Borel subgroup $B$ satisfies 
$(*)$ for $N \ge 2g-1$ then for 
each positive root $\alpha$ the line bundle $(-\alpha)_*(E_{\sigma})$ is 
globally generated and satisfies $H^1(C,\,(-\alpha)_*(E_{\sigma}))=0$.

We now state a basic lemma due to Drinfeld-Simpson \cite{Drinfeld-Simpson}
about the existence of reductions $\sigma$ satisfying the property 
$(*)$ for any $N$.

\begin{lemma}\label{largered} For any $N$ there exists a reduction of 
structure group $\sigma$ of $E$ to $P$ satisfying  the property $(*)$ for $N$
\end{lemma}

\proof The result is shown for the case $P=B$ in Drinfeld-Simpson
\cite{Drinfeld-Simpson} and for the case of arbitrary $P$ one just
observes that giving a reduction of structure group to $B$
automatically gives rise to a reduction of structure group to $P$, by 
extension and
a reduction to $B$ satisfying the property $(*)$ for N will give
rise to a $P$ reduction of $E$ satisfying the same
conditions. $\hfill \square$

\rem We wish to indicate that the proof in \cite{Drinfeld-Simpson}
first reduces the problem to the case when $E$ is a trivial bundle by choosing
a trivialisation on a Zariski open subset of $C$. Then one can reduce
this problem to the case when $G$ is simply connected and $C$ is a
projective line. At this stage one can recover the result by using the
Theorem 7.4 and Proposition 6.13 of Ramanathan
\cite{Ramanathan83}.

Let $S$ be a finite type scheme over $k$. Let $E_L$ be a principal $L$-bundle
over $C\times S$.
There is an conjugation action of $L$ on $R_uP$ using the splitting of 
$P \lra L$. Hence there is a filtration
\begin{equation}
\label{L}
R_uP\,=\, U_1 \, \supset \,U_2 \,\supset\,
\cdots\,\supset\, U_k\,\supset\, U_{k+1}\, = \,\{e\}\, 
\end{equation}
such that $U_j$ is normal in $P$, each quotient $U_j/U_{j+1}$ are invariant
under the action of $L$, lies in the center of $R_uP/U_{j+1}$ and define 
irreducible representations of $L$.

Consider the functor 
$H^1_S(C,\,R_uP(E_L))$ from 
$({\rm Sch}/S)_{{\rm fppf}} \lra {\rm Sets}$ which takes a 
scheme $(s:S'\ra S)$ to isomorphism classes of  pairs $(E_P,\,\phi)$ 
where $E_P$ is a $P$-bundle on $C\times S'$ and $\phi$ is an isomorphism
$\phi : p_*(E_P) \cong s'^*(E_L)$.
We now record the following lemma for later use.
\begin{lemma}\label{unipotent}
The functor $H^1_S(C,\,R_uP(E_L))$ is representable by an affine bundle 
over $S$ under the assumption $H^0(C,\,U_j/U_{j+1}(E_L) |_{C\times\{s\}})=0$
for each $j=1,\ldots,k$ and $s\in S$, where $U_j/U_{j+1}(E_L)$ is a vector
bundle on $C\times S$ associated to the representation of $L$ on 
$U_j/U_{j+1}$.
\end{lemma}
\proof See Theorem A.2.6 of Friedman and Morgan \cite{Friedman-Morgan2} 
for the proof. $\hfill \square$

The following two lemmas hold when the characteristic of the field $k$ is 
zero or when the parabolic $P$ is Borel. The method of proof is similar 
to that of the proof of Lemma 3.6 of Kumar-Narasimhan \cite{Kumar-Narasimhan}.
We will use these results only for the case of Borel subgroups.
\begin{lemma}\label{H0vanish}
Suppose $E$ is a $G$ bundle which admits a reduction of structure group $\sigma$ to $P$ such that the property {\rm ($*$)} holds for $N=1$ 
and that the associated Levi bundle $E_L=p_*(E_{\sigma})$ is semi-stable then 
$H^0(C,\,U_j/U_{j+1}(E_L))=0$ for each $j$.
\end{lemma}
\proof Since the representation of $L$ on $U_j/U_{j+1}$ is irreducible
hence the vector bundle $U_j/U_{j+1}(E_L)$ is
semistable. The degree of this vector bundle is recovered from the
character of $L$ on the highest exterior power of $U_j/U_{j+1}$ which
when restricted to $T$ is non-trivial and is a non-negative linear
combination of simple roots hence by condition ($*$) for $N=1$ we 
have ${\rm deg}(U_j/U_{j+1}(E_L))<0$. This implies that
$H^0(C,\,U_j/U_{j+1}(E_L))=0$.  $\hfill \square$

\begin{lemma}\label{H1vanish}
Let $E$ be a principal $G$-bundle which 
admits a $P$-reduction $\sigma$ such that $[\sigma]$ satisfies the property
{\rm ($*$)} for $N=2g-1$ and that the associated $L$-bundle 
$E_L$ is semistable. Then $H^1(C,\,T_{\sigma})=0$

\end{lemma}
\proof 
Let ${\mf g}$ and ${\mf p}$ be the lie algebras of $G$ and $P$ respectively.
Consider the filtration 
$0=V_0 \subset V_1 \subset \ldots \subset V_k={\mf g}/{\mf p}$ of 
$P$-submodules $V_i$ such that each quotient $V_i/V_{i-1}$ defines a 
irreducible representation of $P$ (hence of $L$). 
Hence we get a filtration of the vector bundle 
$T_{\sigma}$ for the action of $L$ on these successive quotients.
Now again the successive quotients define semistable vector bundles on $C$
and that the character of $T$ defined by restriction of the 
representation of $L$ on $V_i/V_{i-1}$  is a negative linear combination 
of simple roots, hence by stability of $E_L$ and the condition $(*)$ 
for $N=2g-1$ we see that $H^1(C,\, T_{\sigma})=0$. $\hfill \square$

We record here a corollary of the above lemma about smoothness of the space of sections using the deformation theory of the Hilbert schemes.

\begin{corollary}\label{smooth}
If $\E$ is a family of $G$-bundles over a smooth projective curve $\Cu$ 
over a scheme $S$. If $\E$ is a family of $G$-bundles over a smooth projective curve $\Cu$ over a scheme $S$ of genus $g$. Let $[\sigma] \in {\cal X}_*(L)$ be a point satisfying the property {\rm ($*$)} for $N=2g-1$. Let 
$y\in {\rm Sec}^{[\sigma]}_{\E/P}$ be a Hilbert point of a $P$ reduction  
$\sigma$ of $\E _x$ with $x \in S$ such that the 
associated $L$-bundle $p_*E_{\sigma}$ is semistable. Then the natural morphism 
${\rm Sec}^{[\sigma]}_{\E/P}\lra S$ is smooth at $y$.
\end{corollary}

\begin{comment}
Let $E$ be a principal $G$-bundle over $C$ with a reduction $\sigma$ to $P$.
We will denote by ${\rm Gr}_{\sigma}(E)$ the principal $G$ bundle obtained by 
extension of structure group of $E_{\sigma}$ by the composite map 
$P\ra L \hra P \hra G$.
The following proposition is proved in 
Kumar-Narasimhan-Ramanathan \cite{Kumar94} Proposition 3.7(a), when $k$ is 
of characteristic 0 but one can check that it holds in any characteristic.

\begin{proposition}\label{degenerate}
There exists a family of $G$-bundles ${\cal E}$ over $C\times \A ^1$
with the property that ${\cal E}|_{C\times (\sA ^1-0)} \cong {\rm
pr}_1^*(E)$ and ${\cal E}|_{C\times \{0\}}\cong {\rm Gr}_{\sigma}(E)$,
where ${\rm pr}_1$ denotes the first projection $ C\times (\A ^1-0)
\ra C$.
\end{proposition}

\proof 
\end{comment}

Now prove a result which generalizes the Theorem 7.4 of \cite{Ramanathan83}
and has a similar proof.

\begin{proposition}\label{main2}
Let $\sigma$ be a $B$-reduction of a principal $G$-bundle $E$ 
satisfying the property $(*)$ for 
$N=2g$. If $\mu$ be a 1-PS such that $w_0[\sigma]\le \mu$,
then $E$ admits a $P$ reduction $\sigma _0$ with  
$w_0[\sigma _0]= \mu$   
\end{proposition}

The following lemma is an extension of the Lemma 7.4.1 of 
\cite{Ramanathan83} and is a step in the proof of the Proposition \ref{main2}

\begin{lemma}\label{induction}
With the above notations there exists a sequence $w_0[\sigma]=\mu_1$,
$\mu _{2}$, $\ldots$, $\mu_n=\mu$ of elements of ${\cal X}_*(T)$ 
such that  $\mu _{i+1}=\mu_{i}+{\hat \alpha}_{j_i}$ for some 
$\alpha_{j_i} \in \Delta$ with $(\mu_{i},\, \alpha_{j_i})\ge 2g-1$.
\end{lemma}

\proof We set $w_0[\sigma]=\mu_1$. We prove the lemma by a downward
induction.  If for each $l>i$, $\mu _{l}$ has been chosen such that 
$\mu_1 \le \mu _l$ and $\mu _{l+1}=\mu_{l}+{\hat \alpha}_{j_l}$ for some 
$\alpha_{j_l} \in \Delta$ with $(\mu_{l},\, \alpha_{j_l})\ge 2g-1$ 
then we want to make
a choice for $\mu _i$.  Since $\mu_1 \le \mu _{i+1}$ 
we can write $\mu _{i+1}= \mu_1 + \sum k_m {\hat
\alpha}_m$ with ${\alpha}_m \in \Delta$ and $k_m \ge 0$.
Now we want to get rid of one of the ${\hat \alpha}_m$ with $k_m >0$.

First we claim that $ (\sum k_m {\hat\alpha}_m,\, \alpha) >0$
for some $\alpha \in \Delta$ such that $k_m >0$.
If not then we would have $k_m= (\sum
k_m {\hat\alpha}_m,\, w_m) \le 0$ for every dominant weight $w_m$ as the
dominant weights are non-negative rational linear combinations of the
simple roots and the fact that $({\hat \alpha}_m,\, \alpha _n)<0$ for 
$m\neq n$. 
For such an $\alpha$, using the condition $(\mu _1,\,\alpha)\ge 2g$ we get
$(\mu _{i+1},\, \alpha) >2g$, hence 
$(\mu _{i+1}-{\hat \alpha},\, \alpha) \ge 2g-1$ as 
$({\hat \alpha},\,\alpha)=2$. 
Now we define $\mu _i=\mu _{i+1}-{\hat \alpha}$ and the lemma follows by 
induction.
$\hfill \square$

The following lemma along with the preceding one implies the 
Proposition \ref{main2}.

\begin{lemma} \label{1step}
Let $\mu, \,\nu$ be 1-PS such that $\mu = \nu + {\hat \alpha}
$ for some $\alpha \in \Delta$ and $(\nu,\,\alpha)\ge 2g-1$. Let $E$
be a principal $G$-bundle with a reduction of structure group $\sigma$
to $B$ with the property that $w_0[\sigma]=\nu$. Then $E$ admits a
$B$ reduction $\sigma _0$ such that $w_0[\sigma _0]=\mu$
\end{lemma}

\proof Let $P_{\alpha}$ be the minimal parabolic containing $B$
defined by the simple root $\alpha$. Let $P_{\alpha}=
R_{u}(P_{\alpha})\,L$ be its Levi decomposition. Let $Z_{\alpha}=({\rm
ker}\,\alpha)^0$ be the connected component of the kernel of
$\alpha$. Then $Z_{\alpha}$ is the connected component of the center
and ${\ov L}= L/Z_{\alpha}$ is a rank 1 semisimple group. 
We have  the projection $P_{\alpha}\ra {\ov L}$ which induces an isomorphism
$P_{\alpha}/B \cong {\ov L}/{\ov B}=\P^1$. The root
$\alpha$ induces a ${\ov \alpha}\in {\cal X}^*(T/Z_{\alpha})$ which is
the simple root of the ${\ov L}$. The coroot ${\hat {\ov \alpha}}$ on
${\ov L}$ is the image of ${\hat \alpha}$ under ${\cal X}_*(T)\ra
{\cal X}_*(T/Z_{\alpha})$. Similarly we have ${\ov \mu}$ and ${\ov \nu}$
as the images of $\mu$ and $\nu$ under ${\cal X}_*(T)\ra
{\cal X}_*(T/Z_{\alpha})$.

Now for any $\nu \in{\cal X}_*(T)$ the integer $(\nu,\,\alpha)$ is
determined by the composite $\alpha \circ \nu:\G_m\ra \G _m$ which 
takes $z \mapsto z^{(\nu,\,\alpha)}$. Since the map $\alpha :T \lra \G _m$ factors
through $T/Z_{\alpha}$, we see that $(\nu,\,\alpha)=({\ov \nu},\,{\ov
\alpha})$.

Let $E_{\sigma}$ be the $B$-bundle associated to $\sigma$ as in the
statement of the lemma satisfying $w_0[\sigma]=\nu$. Then
$E_{\sigma}$ gives rise to a $P_{\alpha}$-bundle by extension of
structure group which we denote by $E_{\sigma,\alpha}$. We also denote
${\ov E}_{\sigma,\alpha}$ the ${\ov L}$-bundle obtained by the
extension of structure group $P_{\alpha}\ra {\ov L}$.

By Lemma \ref{compare} there is a 
bijective correspondence between the $B$ reductions of
$E_{\sigma,\alpha}$ of numerical type $[\sigma _1]$ and 
${\ov B}$ reductions of ${\ov E}_{\sigma,\alpha}$ of numerical type 
$[{\ov \sigma}_1]$.  Let $\sigma _1$ be a $B$ reduction of
$E_{\sigma,\alpha}$ and ${\ov \sigma}_1$ be the ${\ov B}$ reduction of
${\ov E}_{\sigma,\alpha}$ under the above correspondence.  
Now again by arguments of the lemma \ref{compare} we see that 
 $[\sigma _1]-\nu$ is a multiple of the coroot ${\hat \alpha}$

Hence we are reduced to proving the Lemma
for the case of rank $1$ semi-simple groups.  This case follows from
the following lemma as $PGL(2)$-bundles on curves come from vector
bundles.

\begin{lemma}Let $V$ be a rank two vector bundle which can be written 
as an exact sequence
$$ 0\ra L_1 \ra V \ra L_2\ra 0$$ such that 
${\rm deg}(L_2)\,-\,{\rm deg}(L_1)=m>2g-2$. 
Then $V$ can be written as an exact sequence
 $$ 0 \ra L_1' \ra V \ra L_2'\ra 0$$ such that 
${\rm deg}(L_2)\,-\,{\rm deg}(L_1)=m+2$.
\end{lemma} 
\proof Consider the family ${\cal V}$ of vector bundles on $C\times \A ^1$
such that ${\cal V}_{C\times \{ 0 \} } =L_1\oplus L_2$ and 
${\cal V}_{C \times  \{ x \} }=V$ for $x\in \A ^1-0$.
The condition $m>2g-2$ ensures that ${\rm Sec}^m_{{\cal V}/B} \ra \A ^1$ is 
smooth. Hence it is enough to prove the result for the case 
$V=L_1 \oplus L_2$. In this case we can construct a sub-bundle 
$L_1(-x) \hra L_1\oplus L_2$ by choosing a section
of $L_2 \otimes L_1^{-1}(x)$ which does not vanish at $x \in C$. 
$\hfill \square$

Now the Proposition \ref{main2} follows by the Lemma \ref{induction} and 
Lemma \ref{1step}.$\hfill \square$

\rem \label{field} One observes that the Proposition \ref{main2} holds even if the ground field is not algebraically closed on the condition that $C$ has a rational point.

In the rest of the section we give some applications of the above result.

We fix the above notation. Let $B_0$ be the opposite Borel subgroup
containing $T$ and the negative root spaces. 

In the special case when the Harder-Narasimhan reduction of a principal 
$G$-bundle is
defined on a Borel subgroup satisfying some conditions we answer a question
of Friedman-Morgan (see \cite{Friedman-Morgan1})
on the behavior of the strata with respect to deformation.

\begin{proposition}
Suppose $\mu$ and $\nu$ are two dominant 1-PS satisfying the conditions
$\nu (\alpha)\ge 2g$ for each $\alpha \in \Delta$ and $\nu \le \mu$.
Then there exists a principal $G$-bundle $E$ whose canonical reduction
$\sigma _0$ is a reduction on a Borel subgroup $B_0$ satisfying 
$[\sigma _0]=\nu$ and there exists a family ${\cal E}$ of $G$-bundles on
$C\times \A^1$ such that for ${\cal E}|_{C\times (\sA ^1-0)}\cong {\rm
pr}_1^*(E)$ and ${\cal E}|_{C\times {0}}$ is a principal $G$ bundle
whose canonical reduction $\tau _0$ is again defined on the Borel
subgroup $B_0$ satisfying $[\tau _0]=\mu$
\end{proposition}

\proof Let $\mu$ and $\nu$ be as in the statement of the Proposition.
Let $E_T$ be a principal $T$ bundle such that $w_0[E_T]=\nu$. Let $E$ be
the principal $G$ bundle obtained by extension to $G$. Let $B$ and
$B_0$ be the opposite Borel subgroups intersecting exactly on $T$.  We
then get by extending to $B$ and $B_0$, reductions of structure
groups $\sigma$ and $\sigma _0$ of $E$ to $B$ and $B_0$ respectively
satisfying $w_0[\sigma]=[\sigma _0]=\nu$ and 
$\sigma _0$ being the canonical reduction of $E$ (by the conditions
stated in the proposition). Now we have a $G$-bundle
$E$ with a reduction $\sigma$ to $B$ and a $\mu$ such that 
$w_0[\sigma]=\nu\le \mu$. Hence by Theorem \ref{main2} 
we get a $B$ reduction $\tau$ of $E$ such that 
$w_0[\tau]=\mu$. 
Now by Proposition 3.7(a) of Kumar-Narasimhan-Ramanathan \cite{Kumar94}
we find a family 
${\cal E}$ of $G$-bundles on
$C\times \A^1$ such that for ${\cal E}|_{C\times (\sA ^1-0)}\cong {\rm
pr}_1^*(E)$ and 
${\cal E}|_{C\times {0}}\cong E_1$ where
Now it is easy to see that $E_1$ is a $G$-bundle 
whose canonical reduction $\tau _0$ defined on the Borel
subgroup $B_0$ satisfying $[\tau _0]=\mu$. $\hfill \square$

Next we address the question of characterization of algebraic families
of $G$-bundles on $C$ and irreducibility of the moduli spaces. 
This part is actually mentioned in
Ramanathan \cite{Ramanathan83} without proofs. We also wish to make 
reference to
Drinfeld-Simpson \cite{Drinfeld-Simpson} where the irreducibility of the 
moduli stack
of $G$-bundles are proved. But the basic results here are a little more 
explicit using the boundedness theorem in \cite{Holla-Narasimhan}.

Recall the following definition. Two principal $G$ bundles $E$ and $F$ are 
algebraically equivalent if there is a connected variety $S$ of finite 
type over $k$ and a family ${\cal E}$ of $G$-bundles on $C\times S$ and two 
$k$-valued points $s_0$ and $s_1$ such that ${\cal E}|{C\times {s_0}}\cong E$ 
and ${\cal E}|{C\times {s_1}}\cong F$.

Let ${\hat Q}$ denote the lattice of coroots in ${\cal X}_*(T)$.
When $k=\C$, it is well known that ${\cal X}_*(T)/{\hat Q}$ classifies 
principal $G$-bundles
topologically. What we prove here is the algebraic classification of
principal $G$-bundles in arbitrary characteristic. 

We first state a preliminary lemma which is easy to prove and 
 will be used in the sequel

\begin{lemma}\label{largeroots}
If $\mu _1$ and $\mu_2$ are  1-PS in $T$ such that their images in 
${\cal X}_*(T)/{\hat Q}$ are equal then there is a 1-PS $\mu$ such that
$\mu _i \le \mu$ for $i=1,\,2$ and that $\mu (\alpha)>2g-1$ for each 
$\alpha \in \Phi ^{+}$.
\end{lemma}

Given a principal $G$-bundle $E$ and a reduction $\sigma$ of $E$ to the 
Borel $B$ we define $c(E)_{\sigma}$ to be the image of $[\sigma]$ in
${\cal X}_*(T)/{\hat Q}$.

\begin{lemma}\label{welldefined}
For any two reductions $\sigma$ and $\tau$ of $E$ to $B$, we have 
$c(E)_{\sigma}=c(E)_{\tau}$. In other words $c(E)_{\sigma}$ is independent 
of the reduction $\sigma$.  
\end{lemma}
\proof The idea of the proof is already there in Proposition 6.16 of 
Ramanathan \cite{Ramanathan83}.
To prove the lemma we may assume that $G$ is semi-simple. This is because for any character $\chi \in {\cal X}_*(G)$ we have $\chi _*(E_{\sigma})= 
\chi _*(E_{\tau})=\chi _*(E)$.
Let ${\wt G}$ be the simply connected covering group. Let
$Z$ be the kernel of ${\wt G}\lra G$. Let ${\wt T}$ and 
${\wt B}$ be the maximal torus and the Borel subgroups of ${\wt G}$ which are the inverse images of $T$ and $B$ respectively.
As in the proof of the Proposition 6.16 of \cite{Ramanathan83}, we have a commuting diagram with horizontal rows exact in the flat topology on $C$ 
(and not exact in the category of group schemes)
$$
\begin{array}{ccccccccc}
1 & \lra & Z & \lra & {\wt G} & \lra & G & \lra & 1 \\
  &    &\vequal&  & \uparrow & & \uparrow & & \\  
1 & \lra & Z & \lra & {\wt B} & \lra & B & \lra & 1 \\  
  &    &\vequal&  & \downarrow & & \downarrow & & \\
1 & \lra & Z & \lra & {\wt T} & \lra & T & \lra & 1 

\end{array} 
$$

Then this gives rise to the following commuting diagram of flat cohomologies 
with horizontal rows exact.
$$
\begin{array}{ccccc}
H^1(C,\,{\wt G}) & \lra & H^1(C,\,G) & {\stackrel {\delta _1}{\lra}} &  
H^2(C,\,Z) \\
\uparrow         &      & \uparrow j_* &            & 
\vequal   \\
H^1(C,\,{\wt B}) & \lra & H^1(C,\,B) & \lra &  
H^2(C,\,Z) \\
\downarrow       &      & \downarrow p_* &            & 
\vequal   \\
H^1(C,\,{\wt T}) & {\stackrel {q_*}{\lra}} & H^1(C,\,T) & 
{\stackrel {\delta _3}{\lra}} &  H^2(C,\,Z) \\
\end{array}
$$
Now if $\sigma$ and $\tau$ are two reductions of $E$ to $B$ then the corresponding classes $E_{\sigma}$ and $E_{\tau}$ in $H^1(C,\,B)$ satisfies 
$j_*(E_{\sigma})=j_*(E_{\tau})$ hence from the commutativity of the 
diagram and the fact that $H^1(C,\,T)$ is a group we have
 $\delta_3(p_*(E_{\sigma})-p_*(E_{\tau}))=0$.
The exactness of the bottom row implies existence of an element 
$a \in H^1(C,\, {\wt T})$ such that $q_*(a)=p_*(E_{\sigma})-p_*(E_{\tau})$.
Now it is easy to see from the definitions that the 
following diagram commutes 
$$
\begin{array}{ccc}
H^1(C,\, {\wt T})   & {\stackrel {q_*}{\lra}} & H^1(C,\, T) \\
\downarrow  [\cdot] &                         & \downarrow [\cdot] \\
{\cal X}_*({\wt T}) & {\stackrel {q_*}{\lra}} & {\cal X}_*(T)
\end{array}
$$

Hence we see that $[\sigma]=[p_*(E_{\sigma})]$ and 
$[\tau]=[p_*(E_{\sigma})]$ differ by an element of ${\cal X}_*({\wt T})$. 
Since $q_*({\cal X}_*({\wt T}))$ is exactly the lattice generated by the 
coroots we are through with proof of the lemma.$\hfill \square$

The above lemma enables us to define a map $c$
form the isomorphism classes of principal $G$ bundles to 
${\cal X}_*(T)/{\hat Q}$ by assigning $c(E)=c(E)_{\sigma}$ for any 
reduction $\sigma$ of $E$ to $B$. We call $c(E)$ the topological type of 
$E$.
\begin{proposition} \label{algeq}
Two principal $G$ bundles $E$ and $F$ are algebraically 
equivalent if and only if $c(E)\,=\,c(F)$.
\end{proposition}  

\proof To prove the ``only if'' part of the proposition we need to only 
consider the case when the  
principal $G$-bundles $E$ and $F$ sit in an irreducible family.
Let $S$ be an irreducible finite type scheme over $k$ with two points 
$s_0$ and $s_1$. 
Let ${\cal E}$ be a family of $G$-bundles on $C\times S$ 
such that ${\cal E}|_{C\times {s_0}}\cong E$ and 
${\cal E}|_{C\times {s_1}}\cong F$.
Now choose a reduction of structure of structure group $\sigma$ of $E$
to $B$ satisfying the property $(*)$ for $N=2g-1$.  
Then by Corollary \ref{smooth} the morphism 
${\rm Sec}^{[\sigma]}_{{\cal E}/B} \ra S$ is smooth and the image contains 
the point $s_0$ defined by $E$, hence it contains a neighborhood of $s_0$.
This gives us an open subset $U$ of $S$ with the property that for 
$s \in U$ the bundle ${\cal E}|_ {C\times {S}}$ admits a $B$ reduction of
type $[\sigma]$.  Similarly by choosing a reduction $\tau$ of $F$
satisfying the condition $(*)$ for $N=2g-1$ we get another open
subset $V$ of $S$ such that for each $s' \in V$, the bundle 
${\cal E}|_ {C\times {S}}$ admits a $B$ reduction of type $[\tau]$. 
Since $S$ is irreducible, $U$ and $V$ have non-trivial intersection, 
hence we produce a principal $G$-bundle $E'$ 
which admits a reduction to $B$ with 
numerical types $[\sigma]$ and $[\tau]$. Now by Lemma \ref{welldefined} we 
have
$c(E)\,=\,c(E)_{\sigma}=c(E')_{\sigma}\,=\,c(E')_{\tau}\,=\,c(F)_{\tau}\,
=\,c(F)$. 
This proves the only if part of the proposition.

For the other part of the proposition the main idea is the construction of  
of the family of $B$-bundles which extend to a fixed $T$-bundle.
Let $E$ and $F$ be principal $G$ bundles such the $c(E)=c(F)$. Choose
reductions of structure group $\sigma$ and $\tau$ of $E$ and $F$
respectively satisfying $(*)$ with $N=2g$.  Since $[\sigma]-[\tau]$
is an integral linear combination of coroots, by Lemma
\ref{largeroots} there is a dominant 1-PS $\mu$ with the
property that $w_0[\sigma] \le \mu$ and $w_0[\tau] \le \mu$.  Now by
Theorem \ref{main2} there are reductions $\sigma _0$ and
$\tau _0$ of $E$ and $F$ respectively with the property that 
$w_0[\sigma _0]=w_0[\tau] = \mu$.  Let $S$ be the moduli space of 
$T$-bundles of type $w_0 \mu$. The space $S$ is essentially a product of
Jacobians. There is also a universal family ${\cal E}_T$ over $C\times
S$.  
The conjugation action of $T$ on $U$ using the splitting of 
$B \lra T$ gives us a filtration as in (\ref{L}). The 
condition ($*$) for $N=1$ along with Lemma \ref{unipotent} and 
Lemma \ref{H0vanish} implies that the functor 
$H^1_S(C,\,U({\cal E}_T))$ is representable by an affine bundle 
${\cal H}$ over $S$. Let
${\cal B}$ be the universal family of $B$ bundles over $C\times {\cal
H}$. Now this universal family when extended to $G$ gives a family
of $G$-bundles over a finite type irreducible scheme ${\cal H}$ (in fact smooth)
containing both $E$ and $F$.
This completes the proof of the proposition. $\hfill \square$

In the following Lemma we relate the invariant $c$ with the degree of the $G$-bundle.
Consider the natural map 
$g:{\cal X}_*(T)\lra {\cal X}_*(G):={\rm Hom}({\cal X}^*(G),\,\Z)$ defined by the dual of the restriction map $ {\cal X}^*(G) \ra {\cal X}^*(T)$.
\begin{lemma}\label{degree}
The homomorphism $g$ factors through ${\hat Q}$ to give a homomorphism 
$g: {\cal X}_*(T)/{\hat Q} \ra {\cal X}_*(G)$.
\end{lemma}
\proof 
for the proof we have to show that any coroot ${\hat \alpha}$ with 
$\alpha \in \Phi$ acts trivially on any character 
$\chi$ of $G$. This follows from 
the definition of the coroot as a homomorphism from the maximal torus of 
$SL_2$ to $G$. Now the fact that any character on $SL_2$ is trivial implies 
the result. $\hfill \square$

\rem Recall that the composition $d \circ c(E)$ is equal to the degree 
of the principal $G$-bundle as defined in \cite{Holla-Narasimhan}.

With the above notations Let $P=P_I$ be the parabolic containing the 
Borel subgroup $B$ corresponding to the subset $I \subset \Delta$. 
Let ${\ov L}=L/Z_0(L)$ be the quotient of the Levi by the connected 
component of the center. We denote 
the projection map $P\ra {\ov L}$ by $p$. Let ${\hat Q}_{\ov L}$ be the
coroot lattice of ${\ov L}$. Then ${\hat Q}_{\ov L}$ is the lattice 
generated by
the coroots ${\hat \alpha}$ with $\alpha \in I$. We also denote by ${\ov B}$
(resp. ${\ov T}$) the Borel subgroup (maximal torus) of ${\ov L}$ 
defined  the images of the respective objects under $p$.
We use the notation $c_G$ 
and $c_L$ for the topological type $c$ to indicate difference when we are 
working with $G$-bundles and ${\ov L}$-bundles respectively.

The following is a lemma which is an algebraic version of the Lemma 2.4 of 
Friedman-Morgan \cite{Friedman-Morgan1} whose proof is also very similar.

\begin{lemma}\label{G-L}
If $E$ is a principal $G$-bundle which admits a $P$ reduction 
$\sigma$ then the numerical type $[\sigma]$ and the element 
$c_G(E)$ determines the element $c_L(p_*E_{\sigma})$.
\end{lemma}
\proof
Suppose $E$ and $F$ are two principal $G$-bundles which are
algebraically equivalent. Let $\sigma$ and $\tau$ be $P$-reductions of
$E$ and $F$ respectively.  Let ${\ov \sigma}_0$ and 
${\ov \tau}_0$ be further reductions of the ${\ov L}$-bundles 
$p_*E_{\sigma}$
and $p_*F_{\tau}$. Then these determine reductions $\sigma _0$ and
$\tau _0$ of $E$ and $F$ to $B$ (by Lemma \ref{compare}). 
With the conditions $c_G(E)=c_G(F)$
and $[\sigma]=[\tau]$ we have to show that
$c_L(p_*(E_{\sigma}))=c_L(p_*(F_{\tau}))$.  The condition
$[\sigma]=[\tau]$ implies that the difference $[\sigma _0]-[\tau _0]$
is in the kernel of the natural map 
${\cal X}_*(T)\otimes \Q \lra {\cal X}_*(P)\otimes \Q$ 
which is exactly equal to $ {\hat Q}_L\otimes \Q$. 
Hence $[\sigma _0]-[\tau _0]$ is a rational linear
combination of elements of ${\hat \alpha}$ with $\alpha \in I$. Now the
condition $c_G(E)=c_G(F)$ implies that $[\sigma _0]-[\tau _0]$ is an
integral linear combination of elements ${\hat \alpha}$ with 
$\alpha \in \Delta$. This proves the lemma.$\hfill \square$

\begin{comment}
The above lemma allows us to define $c_L([\sigma]) \in {\cal X}_*({\ov L})$
for a numerical type $[\sigma] \in{\cal X}_*(P)$ for which the 
$G$-bundle $E$ admits a $P$ reduction.
\end{comment}

Recall the definition of the instability degree of a principal $G$-bundle $E$
$$
{\rm Ideg}_G(E)= {\rm Max} \{ {\rm deg}({\rm ad}E_{\sigma}) |
(P,\sigma) \},
$$
where the maximal is taken over all parabolic reductions of the principal 
$G$-bundle $E$. It follows from Behrend \cite{Behrend} that if 
$(P,\sigma)$ is the Harder-Narasimhan reduction of $E$ then 
${\rm Ideg}_G(E)=-d([\sigma])$ 
(also see Mehta-Subramanian \cite{Mehta-Sub} or 
Holla-Biswas \cite{Holla-Biswas} for an account).

Let $c\in {\cal X}_*(T)/{\hat Q}$ be a fixed class. 
Let $M_G(c,d)$ be the set of isomorphism classes of principal $G$ bundles 
over $C$ such that the instability degree is bounded by $d$.   
We need the following proposition.
\begin{proposition}\label{bound}
$M_G(c,d)$ is bounded. In other words there exists a finite scheme $S$ and 
a principal $G$-bundle $\E$ over $S$ such that every member in
$M_G(c,d)$ occurs in $S$.
\end{proposition} 
This is exactly the Theorem 1.2 of \cite{Holla-Narasimhan} 
when $d=0$. For $d >0$ the same method works.

The following is a version of Proposition 4.1, Ramanathan \cite{Ramanathan0}
about openness of stable bundles ($d<0$ case) which is proved in the 
analytic setup but 
the proof goes through in any characteristic.  
\begin{proposition}\label{open}
In any family of bundles $G$-bundles ${\cal E}\ra C\times S$ 
the subset $S^d$ corresponding points $x \in S$ for which 
${\cal E}|_{C\times \{x\}}$ has instability degree less than $d$
is open.
\end{proposition}
The above can also be proved using the Lemma \ref{lowerbound}, the 
properness of the Hilbert schemes and the Lemma \ref{bend-break} 
(to be proved later).

The idea of the proof of the Proposition \ref{algeq} gives something stronger 
which will be used later.
\begin{theorem}\label{irred}
There exists an irreducible smooth variety $S$ and a 
family of $G$-bundles on $C\times S$ such that every element in $M_G(c,d)$ 
occurs in the family $\E$.
\end{theorem}
\proof By Lemma \ref{bound} there is a finite type scheme $X$ with family
of $G$-bundles ${\cal E}$ over $C\times X$ such that for each $x\in X$
the bundle ${\cal E}|_{C\times \{x\}}$ lies in $M_G(c,d)$ and every
$G$-bundle in $M_G(c,d)$ is isomorphic to a member in the family ${\cal E}$.

Consider the defining morphism 
$\bigcup {\rm Sec}^{[\sigma]}_{{\cal E}/B} \lra X$,
where the union is taken over all $[\sigma]$ satisfying the property ($*$)
for $N=2g$. We see that this is a smooth morphism by Corollary 
\ref{smooth}. Since $X$ is of finite type there exists finitely 
many open sets $\{U_i \}_{i=1 \ldots m}$ which cover $X$ and elements 
$[\tau _i]$ for $i=1\ldots m$ such that for each $x \in U_i$ the principal 
$G$-bundle ${\cal E}|_{C\times \{ x \}}$ admits a $B$ reduction with 
numerical type 
$[\tau _i]$ and such that each $[\tau _i]$ satisfies the condition 
($*$) for $N=2g$. 
By Lemma \ref{largeroots} we can find a 1-PS $\mu$ with 
$\mu (\alpha)\ge 2g-1$ for each $\alpha \in \Phi ^{+}$ satisfying 
$w_0[\tau _i] \le \mu$ for each $i$.
Now by Theorem \ref{main2} we see that every  
$G$-bundle which  corresponds to a point in $X$ 
admits a $B$ reduction with numerical type $w_0 \mu $. 

Take $S_0$ to be the moduli space of $T$-bundles of type $w_0\mu$ and
continue with the steps of the ``if'' part of the proposition
\ref{algeq} to obtain an affine bundle $S$ over $S_0$ which 
satisfies the conditions of the statement of the theorem. $\hfill \square$

Using similar arguments and the remark \ref{field}, one can also obtain a 
slight generalization of the Theorem 1 of Drinfeld and Simpson 
\cite{Drinfeld-Simpson}.

\begin{corollary}\label{ds}
If ${\cal E}$ is a principal $G$-bundle over smooth curve 
$\Cu$ over a finite type scheme $S$. Then for any $N$ there is a 
$[\sigma] \in {\cal X}_*(T)$ and a surjective \'etale cover $S' \ra S$ such that 
$[\sigma]$ satisfies the property {\rm ($*$)} for $N$, and $\E$ admits a 
$B$-reduction over $S'$ with numerical type $[\sigma]$.
\end{corollary}

Suppose $P$ is a connected linear algebraic group which is written as 
$R_uP \cdot L$ with the property that $L$ is reductive and $R_uP$ being the 
unipotent radical. We have the natural projection $p: P \ra L$.
We record here a fact which will be used later.

\begin{proposition}\label{versal}
Let $E$ be a principal $P$-bundle over $C$. 
There is an algebraic versal family of deformations of $E$ and the tangent 
space to the deformation functor is $H^1(C,\,{\rm ad}E)$.
\end{proposition}
\proof 
 First we observe that the deformation functor $D_E$ satisfies the 
properties $H_1$, $H_2$ and $H_3$ of Schlessinger \cite{Schlessinger}. 
Hence a formal versal hull exists. 
To verify that the there is a an algebraic versal hull we use the results 
of Artin 
\cite{Artin}. Also the second condition of Artin namely that 
$D_E({\hat A}) \ra \liminv ~ D_E(A/m^n)$ is 
injective for a local ring of an algebraic scheme with residue field $k$, 
is easy to check using a faithful representation of $P$ in $GL_N$ 
and considering the principal $P$-bundle as a vector bundle with a section
on some tensor power of it and using the Grothendieck existence theorems 
\cite{EGA} III 5. This part is done in \cite{Ramanathan83} Proposition 8.4.
So it is enough to check that there is a principal $P$-bundle ${\cal V}$ on 
$C\times S$ for some $S$ smooth scheme of finite type such that the 
infinitesimal deformation map is surjective. 
When the group $R_uP=\{e\}$ that is when $P$ is reductive this follows 
from Proposition \ref{irred} and Corollary \ref{smooth}

In the case $P$ is not reductive we use the filtration (\ref{L}) on
$R_uP$. Now the proposition follows from the following lemma by the
induction on the length of the filtration and by writing the exact
sequence $\{ e \} \ra U_k \ra P \ra P/U_{k} \ra \{ e \}$ and finally
reducing this to the case when the group $P$ is reductive.

\begin{lemma}
Let $e \ra U \ra P \ra M \ra e $ be an exact sequence
of connected algebraic groups such that $U=G_a^n$ and that $P$ action of 
$U$ factors through $M$ to give a representation of $M$.
Let $E$ be a $P$-bundle and let $E_M$ be the associated
$M$-bundle under the surjection $p:P\ra M$. Let ${\cal V}_1\lra C\times S_1$ 
be a family of $M$-bundles parameterizes by a finite type smooth scheme $S_1$
which is versal at $E_M$. Then
there exists a family ${\cal V} \lra C\times S$ of $P$-bundles again 
parameterized by a finite type smooth scheme $S$ which is versal at 
$E$ and there is a surjective morphism $S\ra S_1$.
\end{lemma}
\proof The unipotent group scheme $U({\cal V}_1)$ over $C\times S_1$
actually comes from a vector bundle $W$ over $C\times S_1$. 
If $K^0\ra K^1$ is a complex of vector bundles defined by the 
semi-continuity theorem which computes $R^i({\rm pr}_2)_*W$ then  
one can check that on the vector bundle $K^1$ we have a family of  
$P$-bundles which has the properties mentioned in the lemma.
$\hfill \square$ 

Another fact which we will need is the following.
\begin{proposition}\label{exis}
If the genus of the curve $C$ is atleast two then for any 
$c\in {\cal X}_*(T)/{\hat Q}$ there exists a stable $G$ bundle whose topological type is $c$.
\end{proposition}
\proof

Let $E$ be any principal $G$-bundle over $C$. 
By Proposition \ref{versal} we have a family $\E$ of $G$-bundles 
parameterized by a finite type 
smooth and integral scheme $S$ which is miniversal at every 
point of $S$ and there is an $x\in S$ such that $E=\E_x$. If $E$ admits a 
$P$-reduction $\sigma$ with $d([\sigma])\le 0$ then we have to show that the 
image of the natural map $f: Y={\rm Sec}^{[\sigma]}_{E/P} \ra S$ does not 
contain an open subset of $S$. 

Let $\E _P$ be the family of $P$-bundles on $Y$ obtained by 
the universal property of $Y$ and let $y\in Y$ be such that 
$\E _{P,y}=E_{\sigma}$. 
Let ${\cal V}\lra C\times S_P$ be a 
family of $P$-bundles, parameterized by finite type smooth scheme $S_P$,
which is miniversal at the point $y'=E_{\sigma}\in S_P$
Now by going to an \'etale neighborhood $V$ of $y \in Y$  and of $S_P$,
and an automorphism of $U$ in a neighborhood of $x \in U$ we may assume that 
the restriction of $f$ (again denoted by $f$) defines a 
morphism $f:V \ra S$ which can be written as $j \circ g$, where 
$g:V \ra S_P$ is defined by the versal property of $S_P$ at $E_{\sigma}$ 
and $j:S_P \ra S$ by the versal property of $U$ at $E$.
Hence it is enough to show that
$H^1(C,\,{\rm ad}E_{\sigma}) \ra H^1(C,\,{\rm ad}E)$ is not surjective.
The last statement follows because $H^1(C,\, T_{\sigma})\neq 0$, as 
$d([\sigma]) \le 0$ and $g \ge 2$.   $\hfill \square$     

\section{Hilbert schemes and Moduli space of maps}
In this section we use the compactifications of the space of sections to 
estimate their dimensions. 
For a given reduction $\sigma$ of $E$ to $P$ 
we consider the Hilbert scheme ${\rm Hilb}^{[\sigma]}_{E/P}$ 
and a the open subscheme ${\rm Sec}^{[\sigma]}_{E/P}$ as defined before. 
The following theorem was proved in
Holla-Narasimhan \cite{Holla-Narasimhan} about the dimension estimates 
of the Hilbert schemes corresponding to a minimal section $\sigma$.

\begin{theorem}
Let $\sigma$ be a minimal section. Let $X$ be a irreducible component
of the Hilbert scheme ${\rm Hilb}^{[\sigma]}_{E/P}$ which contains $\sigma$
as a Hilbert point. Then every point in $X$ corresponds to a Hilbert
point of a section. In other words $X \subset {\rm Sec}^{[\sigma]}_{E/P}$.
Moreover ${\rm dim}(X) \le {\rm dim}(G/P)$ and 
${\rm deg}(T_{\sigma})\le g\,{\rm dim}(G/P) $.
\end{theorem}

\rem \label{HN} The dimension bound 
${\rm dim}(X) \le {\rm dim}(G/P)$ follows from the first assertion of 
the above theorem by a rigidity argument which implies that the 
evaluation morphism $X \times C \ra E/P$ is finite. 
(see Lemma 2.4 of \cite{Holla-Narasimhan}). 
The last assertion follows from above by a deformation theoretic argument 
(see Proposition 3 of Mori \cite{Mori}).

In our proof of the fact that the minimal section satisfies ${\rm
deg}(T_{\sigma}) \le g \, {\rm dim}(G/P)$, we only needed the fact that
the highest exterior power of $T_{\pi}$ is ample over the fibers of
the map $\pi : E/P \lra X$. More generally the proof of the above theorem 
would go through assuming $\sigma$ to be numerically minimal.
This is the content of the proposition below.

\begin{proposition}\label{Ex}
In the above set up, if $X$ is an irreducible component of 
${\rm Hilb}^{[\sigma]}_{E/P}$ which contains a Hilbert point of an 
numerically minimal section
$\sigma$ then every element in $X$ corresponds to a Hilbert point of a section.
In other words $X\subset {\rm Sec}^{[\sigma]}_{E/P}$. 
Moreover we have ${\rm dim}(X) \le {\rm dim}(G/P) $ and ${\rm
deg}(T_{\sigma}) \le g \, {\rm dim}(G/P)$.
\end{proposition}

\proof 
We need the following lemma for the proof.
\begin{lemma} \label{bend-break}
If the Hilbert point of a closed subscheme $Y$ of 
$E/P$ is in ${\rm Hilb}^{[\sigma]}_{E/P}$ then 
$Y$ has a unique irreducible component 
$C_0$ which maps isomorphically onto $C$ under the composition 
$C_0 \lra E/P \lra C$, hence defines a section $\sigma _0$ of $\pi$ which 
satisfies the inequality $[\sigma] \le [\sigma _0]$. 
More over if the above is an equality then the subscheme $Y$ coincides with 
$C_0$. 
\end{lemma}
\proof of the lemma:
We just follow the arguments of the Proposition 2.3 of 
\cite{Holla-Narasimhan}. By choosing a line 
bundle $L$ on $C$ of degree one we see that 
$${\chi}\,(Y, { \cal O}_Y) = {\chi}\,(C, { \cal O}_C) ~~~\mbox{and}~~~ 
{\chi}\,(Y, f^*\,(L)) = {\chi}\,(C,L),$$

By applying Lemma 2.2 (i) of \cite{Holla-Narasimhan} we get the unique 
irreducible component $C_0$ mapping isomorphically onto $C$.
The third part of the Lemma 2.2 (iii) now implies  that for any line bundle 
$\xi$ on 
$E/P$ which is ample along the fibers of $\pi$ we have 
$\mbox{deg}\,(C_0,\xi) \le \mbox{deg}\,(Y,\xi)$
and this is an equality if and only if there are no other one
dimensional components. This proves  the Lemma \ref{bend-break} by 
taking the line
bundle $\xi$ to associated to anti-dominant characters of $P$. The final
remark is that the zero dimensional components automatically disappear
once there are no other one dimensional components by Lemma 2.2 (ii)
of \cite{Holla-Narasimhan}. $\hfill \square$

The above lemma immediately implies the first part of the Proposition
\ref{Ex} and the other parts follow form Remark \ref{HN}. $\hfill \square$

\rem \label{proper}
Let $[\sigma]$ be a numerical type not necessarily minimal. 
Suppose that we have an irreducible component $X$ of  
${\rm Hilb}^{[\sigma]}_{E/P}$ which is contained in 
${\rm Sec}^{[\sigma]}_{E/P}$,
in other words every element of $X$ is a Hilbert point of a section.  
Then we see again by Remark \ref{HN} that  
${\rm dim}(X)\le {\rm dim}(G/P)$ and $d([\sigma])\le g\cdot {\rm dim}(G/P)$.
In some sense the minimality condition for the section $\sigma$ is used only 
to get components of Hilbert schemes which do not have any boundaries.

\begin{corollary}\label{finite}
There are only finitely many numerically minimal points in ${\cal X}_*(P)$ 
corresponding to $P$ reductions of the principal $G$-bundle $E$. 
\end{corollary}
\proof This follows from Lemma \ref{lowerbound} and the last inequality 
in the Proposition \ref{Ex}.$\hfill \square$
 
We have the following result on the dimension estimates for the 
irreducible components of the Hilbert scheme containing the Hilbert point 
of a section.

\begin{theorem}\label{hilb}
Let $\gamma _1,\, \gamma_2,\, \ldots , \gamma_m$ be the numerically minimal
types in ${\cal X}_*(L)$ for the principal $G$-bundle $E$.
If $X$ is an irreducible component of ${\rm Hilb}^{[\sigma]}_{E/P}$ which 
contains the reduction of structure group $\sigma$ as a Hilbert point
then there exists a $i$ with $[\sigma] \le \gamma_i$ such that
$${\rm dim}(X) \le {\rm dim}(G/P)+
d({[\sigma]})-d({\gamma_i}).$$
\end{theorem}

The above result is proved using a similar result on the dimension estimates 
of the moduli space of stable maps. Now we recall the basic facts about 
the moduli space of stable maps.

\begin{comment}
The above theorem is a  generalization of the result of 
Mihnea Popa and Mike Roth \cite{Popa-Roth}
on the dimension estimates of Quot scheme 
and its proof follows the same ideas as theirs 
but one has to take care of the numerical types of the sections.
As in their case of Quot schemes, the above theorem follows from a similar 
estimates for the dimension of the moduli space of maps.
Now we recall the basic facts about moduli space of stable maps.
\end{comment}

Let $X$ be a smooth projective variety.  We consider the
Kontsevich space of stable maps ${\ov {\cal M}}_{g,n}(X,
\beta_d)$. As a functor, its $S$-valued points parameterizes
stable families over $S$ of maps from $n$-pointed,
genus $g$ curves to $X$ representing the class $\beta$ with isomorphisms.  
More precisely equivalence class of triples $( \pi : {\cal C}\lra S,\,
\{p_i\}_{1\le i\le n},\,\mu : {\cal C} \lra X)$ satisfying
the following properties.
\begin{enumerate}
\item A family of $n$-pointed, genus $g$ curves $ \pi : {\cal C}\lra
S$ which is flat and projective.
\item $n$ sections $\{p_1,\,\ldots\,,\, p_n\}$ of $ \pi : {\cal C}\lra
S$ such that each geometric fiber \\ $({\cal C}_s, p_1(s),\, \ldots
,\, p_n(s))$ is an $n$-pointed genus $g$ curve which is projective,
connected, reduced, nodal curve of arithmetic genus $g$ with $n$
distinct, nonsingular, marked points.
\item For each geometric point $s\in S$ the morphism 
$\mu :{\cal C} \lra X$ restricted to the fiber $\Cu _s$ 
satisfies the following.
\begin{enumerate}
\item If a rational component $F$ of ${\cal C}_s$ is mapped to a point, then 
$F$ must contain at least three special points (marked points or nodes).
\item If a component $F$ of arithmetic genus 1 is mapped to a point, 
then it must contain at least one special point. 
\item $\beta = (\mu |_{ {\cal C}_s})_*[ {\cal C}_s]$
\end{enumerate}
\end{enumerate}

Here the last equality holds in the homology groups ${\rm H}_2(X,K)$, 
where $K$ is a characteristic zero field (mostly $\Q$ or $\Q _l$). 
If the variety lives in characteristic  $p$ then we replace the above 
homology by the \'etale cohomology group ${\rm H}^{2n-2}(X,\Q _l)$ (or equivalently $H^2(C,\,\Q _l)^* $) for a prime $l$ different from $p$.

${\ov {\cal M}}_{g,n}(X, \beta_d)$ is known to have a structure of a proper
 Artin algebraic stack, with finite automorphism at $k$-valued points, 
which admits a projective coarse moduli space 
$ {\ov  M}_{g,n}(X, \beta_d)$ in all characteristics. For proofs see 
Fulton-Pandharipande \cite{Fulton-Pandharipande}, 
Harris-Morrison\cite{Harris-Morrison}, and 
mostly Abramovich-Oort \cite{Abramovich-Oort}.

We also have an open substack of above giving rise to a coarse moduli scheme 
$M_{g,n}(X, \beta_d)$, which parameterizes  maps from smooth curves.
Furthermore there is a ``forgetful'' map  
${\ov  M}_{g,n}(X, \beta_d)\lra {\ov  M}_{g,n-1}(X, \beta_d)$ which 
extends the natural forgetful map of the open moduli spaces, with one 
dimensional fibers. 

We will be only interested in the case when the space $X$ is of the
form $E/P$ where $E$ is a principal $G$ bundle on a smooth projective
curve $C$, and $P$ is a parabolic subgroup of $G$. Also in our case
we will mostly assume $g=g_C$ and $\beta_d$ is a class in  
$H^2(X,\, K)^*$ determined by a reduction of structure group $\sigma$
of $E$ to $P$.

The basic relation between the homology classes $\beta$ and the 
numerical type $\sigma$ can be described as follows. 

Given a $\beta \in H^2(E/P,\, K)^*$ we define a point 
$[\beta]\in {\cal X}_*(P)$ by 
$[\beta](\chi)=\beta(c_1(L_{\chi}))$ where $L_{\chi}$ is the line bundle 
on $E/P$ defined by $\chi$. This defines a homomorphism 
$H^2(E/P,\, K)^*\ra {\cal X}_*(P)\otimes K$. 

Since $H^2(E/P,\, K)$ is generated by the first Chern classes 
of the line bundles (as $H^2(E/P,\,{\cal O}_{E/P})=0$) and the fact that 
every line bundle on $E/P$ is uniquely of the form $L_{\chi}\otimes \pi ^*(L)$ 
 with a line bundles $L$ over $C$ and $L_{\chi}$ over $E/P$ defined by a character $\chi$ of $P$. Hence we have a well defined set theoretic splitting of the 
above homomorphism defined by $\beta _{[\sigma]}=[\sigma](\chi)+{\rm deg}(L)$ 

When $\sigma$ is a reduction of structure group of a principal $G$-bundle
$E$ then $\sigma$ defines a point of both the spaces 
$M_g(E/P, \beta _{[\sigma]})$ and ${\rm Sec}^{[\sigma]}_{E/P}$.
In the following lemma we prove that the above correspondence defines an isomorphism between the two. 
\begin{lemma} \label{compare2}
There is a natural isomorphism between $M_g(E/P, \beta _{[\sigma]})$ and
${\rm Sec}^{[\sigma]}_{E/P}$ which takes a point corresponding to
the irreducible curve to the section defined by it.
\end{lemma}
\proof We first construct the map from the algebraic stack ${\cal
 M}_g(E/P, \beta _{[\sigma]})\lra{\rm Sec}^{[\sigma]}_{E/P}$, this will
 automatically define the map from the coarse moduli space.  
Given a tuple $(q:{\cal C}\ra S,\,f:{\cal C}\lra E/P) $
 with ${\cal C}$ being a flat family of smooth curves over $S$,
we see that $(q,\,\pi \circ f):{\cal C}\ra S\times C$ 
defines an isomorphism over $S$. Hence this setup  
automatically gives rise to flat family
 of sections. This defines the morphism.  Notice that all the
 isomorphisms in ${\cal M}_g(E/P, \beta_{[\sigma]})$ are automatically
 collapsed in ${\rm Sec}^{[\sigma]}_{E/P}$. The inverse of the above
 morphism can be obtained by simply inverting the above operation and 
then composing it with the natural morphism ${\cal M}_g(E/P, \beta_{[\sigma]})
 \lra M_g(E/P, \beta_{[\sigma]})$. Hence the lemma is proved. 
$\hfill \square$

\rem \label{propermaps} 
Suppose there is a component $X$ of ${\ov M}_g(E/P, \beta_{[\sigma]})$
which has no points of the boundary, in other words, every point of
$X$ corresponds to a map from an irreducible curve. Then 
such a component will map isomorphically on to an irreducible
component of 
${\rm Sec}^{[\sigma]}_{E/P} \subset {\rm Hilb}^{[\sigma]}_{E/P}$. Since $X$
is proper we see that $X$ defines an irreducible component of 
${\rm Hilb}^{[\sigma]}_{E/P}$.
Hence by Remark \ref{proper},
we obtain ${\rm dim}(X) \le {\rm dim}(G/P)$ 
and $d([\sigma])\le g\cdot {\rm dim}(G/P)$.

Now we state the main result of this section
\begin{theorem}\label{stable}Let $\gamma _1,\, \gamma_2,\, \ldots , 
\gamma_m$ be the numerically minimal types in ${\cal X}_*(P)$ for $E$.
Let $[\sigma]$ be a numerical type.
If $X$ is an irreducible component of ${\ov M}_{g}(E/P, \beta_{[\sigma]})$
then there exists an $i$ with $[\sigma] \le \gamma _i$ such that 
$${\rm dim}(X) \le {\rm dim}(G/P)+
d([\sigma])-d(\gamma_i).$$
\end{theorem}

The above result is a generalization of a result of 
Mihnea Popa and Mike Roth \cite{Popa-Roth}
on the dimension estimates of the moduli space of stable maps in the case 
of $G=GL_n$ 
and its proof follows similar ideas 
but one has to take care of the numerical types of the sections.

Suppose $X$ is a component of ${\ov M}_{g}(E/P, \beta_{[\sigma]})$ whose
generic point corresponds to a map with irreducible domain. Let $Y$ be
the boundary in $X$ corresponding to maps with reducible domain. 
Let $Y'$ be an irreducible component of $Y$. 

The following basic lemma which is a stronger analogue of 
Lemma \ref{bend-break} for the moduli space of maps and which 
does not hold for the case of Hilbert schemes.

\begin{lemma}\label{smoothen} A generic
element of $Y'$ corresponds to a map from a connected reduced nodal
curve $C'$ with irreducible components $C_0,\, C_1,\, \ldots,\, C_k$
such that
\begin{enumerate}
\item $C_0$ maps isomorphically onto $C$, every other $C_i$ is isomorphic to 
$\P^1$. 
\item Only singularities of $C'$ are ordinary nodes with each $C_i$, for 
$i=1\ldots k$, intersecting $C_0$ at a point $x_i$ with $x_i\neq x_j$ for 
$i\neq j$. And these are only intersections between the irreducible components.
\item
$C_0$ gives rise to section  $\sigma_0$ of $\pi$ whose numerical type 
$[\sigma _0]$ satisfies $[\sigma] \le [\sigma _0]$. Moreover
for any line bundle $\xi$ on $E/P$ ample along the fibers of $\pi$, we have 
$${\rm deg}(C_0, \xi|_{C_0})+\, \sum_{i=1}^k {\rm deg}(C_i,\xi|_{C_i})=
{\rm deg}(C', \xi|_{C'}),$$
with $ {\rm deg}(C_i,\xi|_{C_i}) >0$ for $i=1\ldots k$.
\end{enumerate}
\end{lemma} 

\proof 
The first part of the proof of this lemma is same as the proof of the Lemma
\ref{bend-break}, but we have to take care of rational tails. 
Since $C'$ is a connected reduced curve with
only nodal singularities, the arithmetic genus of $C'$ coincides with
that of $C$ and $C_0$. This forces the other components to be isomorphic to 
$\P ^1$. Moreover we see that these rational components form trees 
living at the fibers of the morphism $\pi$. 
To show that such a generic curve $C'$ is a comb we 
need to smoothen the rational curves that we
encounter. This follows
from a general statement  proved in Theorem 7.6 II, p. 155 of Kollar 
\cite{kollar}, which uses the fact that $f^*T _{\pi}$ when restricted to 
each of the components of the tails is semi-positive.  This is always the 
case as the tangent bundles of the flag varieties are globally generated and
the tails lie in the fibers of the morphism $\pi$.  $\hfill \square$

\rem The Lemma \ref{smoothen} also hold in 
the situation where $Y=Y'$ is an irreducible component of 
${\ov M}_{g}(E/P, \beta_{[\sigma]})$ consisting of only reducible curves. 
This will also be used in the proof of the Theorem \ref{stable}.

One of the steps in the proof of Theorem \ref{stable} is the following 
lemma which is essentially follows from deformation theory of the nodal curves.

\begin{lemma} \label{codimension} Let $X$ be an irreducible component of 
${\ov M}_{g}(E/P, \beta_{[\sigma]})$ whose general element corresponds to a
irreducible curve. Let $Y' \subset X$ be an irreducible component of
the boundary $Y$ corresponding to the reducible curves. Then
the codimension of $Y'$ in $X$ is at most $k$, where $k$ is the number
of nodes in a curve corresponding to a generic element of $Y'$.
\end{lemma}

\proof This is a standard statement about deformation theory. It is
enough to show that in the algebraic miniversal deformation of such a
curve, the boundary is of codimension $l$. The construction of the
deformation is done in Vistoli \cite{Vistoli}

The tangent space to the deformation  functor is 
${\rm Def}_1(C')= {\rm Ext}^1_{{\cal O}_{C'}}(\Omega_{C'},\,{\cal O}_{C'})$ 
and it is calculated by the following exact sequence.
$$
0\ra H^1(C',\,\theta_{C'}) \ra 
{\rm Ext}^1_{{\cal O}_{C'}}(\Omega_{C'},\,{\cal O}_{C'}) \ra
H^0(C',\, {\underline {\rm Ext}}^1_{{\cal O}_{C'}}
(\Omega_{C'},\,{\cal O}_{C'})) \ra
H^2(C',\,\theta_{C'})
$$
In our situation we have $H^2(C',\,\theta_{C'})=0$. An easy
computation using the normalization of $C'$ shows that ${\rm
dim}\,H^1(C',\,\theta_{C'})=3\,g_{C'}\,-\,3\,-\,k$ and ${\underline {\rm
Ext}}^1_{{\cal O}_{C'}}(\Omega_{C'},\,{\cal O}_{C'})$ is sky scraper
sheaf supported on the nodes and has length $1$ at each of the nodes
of $C'$.  These facts imply that ${\rm Def}_1(C')$ has dimension
$3\,g_{C'}\,-\,3$ containing the subspace ${\rm
dim}\,H^1(C',\,\theta_{C'})$ which parameterizes the deformations of the
curve $C'$ preserving the singularities. 
Note that in the case these dimensions are negative one has to put 
additional markings to ensure that the nodal curves have only finite 
automorphisms and then work over ${\ov M}_{g,n}(E/P,\,\beta_{[\sigma]})$ 
instead of ${\ov M}_{g}(E/P,\,\beta_{[\sigma]})$.
Using this the Lemma follows.  $\hfill \square$

\proof of the Theorem \ref{stable}: 
We first prove the result when the generic element in $X$ 
corresponds to a map from an irreducible domain.
Let $Y'\subset Y$ be an
irreducible component of the boundary $Y$ in $X$. If $(C',\,f)$ is a
curve which corresponds to a general element of $Y'$, then by Lemma
\ref{smoothen}, we can write  $C'= \bigcup _{1=0}^k C_i$ 
where $C_0$ defines a section
$\sigma _0$ of $\pi$ and the other $C_i$'s are isomorphic to $\P^1$.  
Hence the curve  $C_0$ along with the points $\{x_i|i=1\ldots k\}$ defines 
a point in $M_{g,k}(E/P, \beta _{[\sigma _0]})$ and the curves $C_i$ along 
with the point $x_i$ defines an element of $M_{0,1}(G/P,\alpha _i)$,
where $\alpha _i$ is the second homology class defined by
$(f|_{C_i})_*[C_i]$.  We then can estimate the dimension of such curves
by first taking irreducible curves of type $[\sigma _0]$ and $k$
copies of rational curves each with types $\alpha_i$ for each $i$ 
such that sum $[\sigma _0]+\sum_{i=1}^k[\alpha_i]=[\sigma]$.

The dimension of ${\ov M}_{0,1}(G/P, \alpha _i)$ can be computed by
deformation theory to be exactly equal to 
${\rm dim}(G/P)+d([\alpha _i])-2$ (see
Theorem 2, Section 1, \cite{Fulton-Pandharipande}).  
From here it follow that ${\ov M}_{0,1}(E/P, \alpha _i)$ has dimension 
equal to ${\rm dim}(G/P)+d([\alpha _i])-1$.  Now the moduli space 
${\ov M}_{g,k}(E/P,\beta _{[\sigma _0]})$ and 
${\ov M}_{0,1}(E/P, \alpha _i)$ have natural
evaluation maps $ev_i:{\ov M}_{g,k}(E/P,\beta _{[\sigma _0]}) \lra E/P$ and
$ev'_i:{\ov M}_{0,1}(G/P, \alpha _i) \lra E/P$ for each $i=1\ldots k$. 

Consider the space
$$Z([\sigma _0], \alpha_1,\ldots, \alpha _k) =
{\ov M}_{g,k}(E/P, \beta _{[\sigma _0]})\times_{\prod _{i=1}^k E/P}
\prod _{i=1}^k {\ov M}_{0,1}(E/P, \alpha _i).$$

This space has dimension 
$$
\begin{array}{l}
{\rm dim}({\ov M}_{g}(E/P,\, \beta _{[\sigma _0]}))+k+ 
\sum_{i=1}^k{\rm dim}({\ov M}_{0,1}(G/P, \alpha _i))-
k\,{\rm dim}(E/P) \\
={\rm dim}({\ov M}_{g}(E/P,\,\beta _{[\sigma _0]}))+k+k\,{\rm dim}(G/P)+
\sum_{i=1}^kd([\alpha _i])-k-k\,{\rm dim}(G/P)-k\\
={\rm dim}({\ov M}_{g}(E/P, \beta _{[\sigma_0]}))
+\sum_{i=1}^kd([\alpha _i])-k\\
={\rm dim}({\ov M}_{g}(E/P, \beta _{[\sigma_0]}))+ 
d([\sigma])-d([\sigma _0])-k
\end{array}.
$$
The dimension of the above will be related to the dimension of $Y'$
using the natural morphism $h_{([\sigma _0], \alpha_1,\ldots, \alpha
_k)}:Z([\sigma _0], \alpha_1,\ldots, \alpha _k)\lra 
{\ov M}_{g}(E/P, \beta_{[\sigma]}) $ 
which is defined by $\{C_i\} \mapsto C'$.
The existence of such a morphism is easy to see at the level of 
Algebraic stacks by taking the valued groupoids but such a morphism 
would automatically descent to the coarse moduli spaces.

Now by Lemma \ref{lowerbound} and the fact that there are only
finitely many possibilities for the collection
$\{[\alpha_1],\ldots,[\alpha_k] \}$ as they are squeezed between 
$[\sigma]$ and $[\sigma _0]$ and generically every element of $Y'$ is
in the image of one such map hence we obtain an irreducible component $X_0$ of 
some ${\ov M}_{g}(E/P, \beta _{[\sigma_0]})$ and  a collection 
$\{\alpha_1,\ldots,\alpha_k \}$ of elements in $H^2(E/P,\, K)^*$ 
such that the image of the morphism  
$h_{([\sigma _0], \alpha_1,\ldots, \alpha_k)}$ contains an open subset of 
$Y'$. Hence we have  
$$
{\rm dim}(Y')\le {\rm dim}(X_0)+ d([\sigma])-
d([\sigma _0])-k. 
$$
Also the above construction implies that the generic point of 
$X_0$ corresponds to map from an irreducible domain.
 Now by Lemma \ref{codimension} we see that the dimension of $X$ is bounded by
\begin{equation}
{\rm dim}(X)\le {\rm dim}(X_0)+ d([\sigma])-
d([\sigma _0]). \label{inequality}
\end{equation}
The above is the main step for the induction in the proof of the Theorem.
When the induction stops we have two possible cases to take care.

So starting with the numerical type $[\sigma]$ we produce an numerical type
$[\sigma _0]$ and an irreducible component $X_0$ of 
${\ov M}_{g}(E/P, \beta _{[\sigma_0]})$ such that the above inequalities hold.
Now by induction we continue this process to get a sequence
of numerical types $[\sigma _{1}],\,\cdots,\, [\sigma _{l}],\,\cdot \cdot$ 
and irreducible components $X_i$ of 
${\ov M}_{g}(E/P, \beta _{[\sigma_i]})$ for each $i$ such that 
\begin{equation}
{\rm dim}(X_i)\le {\rm dim}(X_{i+1})+ d([\sigma _i])-
d([\sigma _{i+1}]) \label{inequalityi}
\end{equation}
This process goes on  until we reach an $l$ for which 
$X_l$ has no boundary or $[\sigma _l]$ is numerically minimal 
(hence no boundary points).
Now the problem reduces to estimating the dimensions
of the irreducible components which contain no boundary points.  In
this case, by Remark \ref{propermaps}, we have 
${\rm dim} X_l \le {\rm dim}(G/P)$.
Hence combining this with
(\ref{inequalityi}) and (\ref{inequality}) we obtain
$$
{\rm dim}(X)\le {\rm dim}(G/P)+ d([\sigma])-
d([\sigma _l]). 
$$
Now if the component $X$ of ${\ov M}_{g}(E/P, \beta_{[\sigma]})$ consists
completely of boundary points, even then we can carry out the above
induction by taking $Y'=X$ and we would get a sequence of numerical types
 $[\sigma _{0}]\ldots [\sigma _{m}]$ and irreducible components $X_i$ of 
${\ov M}_{g}(E/P, \beta _{[\sigma_i]})$ such that
\begin{equation}\label{ineqbound}
{\rm dim}(X)\le {\rm dim}(X_0)+ d([\sigma])-
d([\sigma _0])-k, 
\end{equation}
and the inequality  (\ref{inequalityi}) holds for $i=0\cdots m-1$. 
Here $k$ is the number of nodes in a curve corresponding to a generic point 
of $X$.
Hence the above argument gives
$$
{\rm dim}(X)\le {\rm dim}(G/P)+ d([\sigma])-
d([\sigma _m])-k, 
$$
which is stronger than the bound for an irreducible component whose generic 
point to a map with irreducible domain.
The Theorem \ref{stable} now follows because  $d([\sigma])-
d([\sigma _l])\le d([\sigma])-d({\gamma _j})$ for some numerically minimal  
point $\gamma_j$.  $\hfill \square$

\proof of the Theorem \ref{hilb} : Let $X$ be an irreducible component of 
${\rm Hilb}^{[\sigma]}_{E/P}$ which contains the Hilbert point of a section 
$\sigma$. Let $X'$ be the open subset of $X$ defined by
$X'=X \cap {\rm Sec}^{[\sigma]}_{E/P}$. Then 
$X'$ is in the image of the isomorphism defined in the 
Lemma \ref{compare2}. Hence Theorem \ref{hilb} follows from the Theorem 
\ref{stable}.

\rem \label{veryuseful}
The proof of the Theorem \ref{stable} shows that if $X$ is an irreducible 
component of ${\ov M}_{g}(E/P, \beta_{[\sigma]})$  for some $[\sigma]$ then 
there exists a $P$ reduction $\sigma _0$ of $E$ 
with $[\sigma _0] \ge [\sigma]$ and  
$d([\sigma _0])\le g \cdot {\rm dim}(G/P)$ 
 such that 
$$
{\rm dim}(X)\le {\rm dim}(X_0)+d([\sigma ])-d([\sigma _0])-k.
$$
Here $k$ is the number of nodes in the curve corresponding to general point of $X$ and $X_0$ is an irreducible component of  
${\ov M}_{g}(E/P, \beta_{[\sigma _0]})$ containing $\sigma _0$. 
This remark will be useful later.

As in the case of \cite{Popa-Roth} we have a 
dimension estimate for the lower bound and this will be used later.

\begin{proposition}\label{lbd}
Let $[\sigma]$ and $[\tau]$ be two numerical types in ${\cal X}_*(L)$ 
such that $[\sigma] \le [\tau]$. Then we have
the following inequality 
$$
{\ov M}_{g}(E/P, \beta_{[\sigma]})\ge {\ov M}_{g}(E/P, \beta_{[\tau]})
+d([\sigma])-d([\tau])-1 
$$
\end{proposition}
\proof 
Let $\alpha =[\sigma]-[\tau]$. Firstly $\alpha$ defines a canonical 
element of $H^2(E/P,K)^*$ (independent of the chosen line bundle $L$ over 
$C$). Now one can check from the definitions that the   
 $h_{([\tau], \alpha)}:Z([\tau], \alpha)\ra 
{\ov M}_{g}(E/P, \beta_{[\sigma]}) $, defined in the proof of the Theorem 
\ref{stable}, is generically injective. Hence the proposition follows from 
the dimension estimate for $ Z([\tau], \alpha)$ as obtained in Theorem 
\ref{stable}. $\hfill \square$

\rem \label{useful} One observes that if $[\sigma]$ is such that 
${\ov M}_{g}(E/P, \beta_{[\sigma]})$ is irreducible and that the generic 
element corresponds to map from an irreducible curve then the proof of the 
Proposition \ref{lbd} actually shows that 
${\ov M}_{g}(E/P, \beta_{[\sigma]})\ge {\ov M}_{g}(E/P, \beta_{[\tau]})
+d([\sigma])-d([\tau]) $, as the image of $h_{([\tau], \alpha)}$ will not 
be generically surjective. This remark will be used later.

\section{Irreducibility and Generic Smoothness of reductions}
We make a temporary change in the notations.
Let $C$ be a smooth projective curve over an algebraically closed field $k$.
Let $G_k$ be a semisimple simply connected algebraic group over $k$.
Let $E$ be a principal $G_k$-bundle over $C$. In this section we want to prove 
that when the numerical type is large enough the space of reductions of 
$E$ to a fixed parabolic subgroup is irreducible and generically smooth.
Our method of proof is to derive this from a similar result for the case when 
the parabolic is a Borel subgroup, and use the results in the previous 
sections to prove it for the parabolic case. The result for the Borel 
subgroups was proved by Harder \cite{Harder} for the case when the curve 
is defined over the algebraic closure of a finite field. From here it does 
not directly follow for arbitrary fields but the method of proof 
works. This will be the first part of this section.

Let $B_k$ be a Borel subgroup of $G_k$. Let $T_k\subset B_k$ be a maximal 
torus. The first result we want to prove in this section is the following.

\begin{theorem}\label{borel} There exists an integer $N$ such that 
if $E$ admits a $B$ reduction of numerical type $[\sigma]$ satisfying the 
property $(*)$ for $N$ then 
${\rm Sec}^{[\sigma]}_{E/B}$ is an irreducible smooth variety of 
dimension
$d([\sigma])+(1-g){\rm dim}(G/B)$. 
\end{theorem}
\proof
Our method of proof is to reduce the problem to the case where we can apply 
the methods of Harder.
We make some first reductions.

We will denote by $F$ the prime field $\F _p$ or $\Z$ depending on whether
our curve is in characteristic $p$ or $0$.
Since the curve $C$ and the principal $G$-bundle is defined by finitely many 
equation, we may assume that there is a finite type affine integral scheme 
$S$ over $F$, a semi simple simply connected split algebraic 
group scheme $G=G_S$ over $S$, a curve ${\cal C}\ra S$ which is smooth and 
proper over $S$ with geometrically integral fibers, 
and there is a principal 
$G$-bundle ${\cal E}$ over ${\cal C}$ such that over the generic 
point ${\rm Spec K}\ra S$, the principal $G_K$-bundle $E_K$ over the 
curve $C_K$ extends to the $G_k$-bundle $E$ over $C$ via the 
field extension ${\rm Spec}(k)\ra {\rm Spec}(K)$.

Recall the definition of the instability degree ${\rm Ideg}_G(E)$ of a 
principal $G$-bundle $E$
We now prove some basic lemmas needed in the proof of the theorem.
\begin{lemma}\label{can}
Let $E$ be a principal $G$-bundle over a smooth projective geometrically connected curve  $C$ over a perfect field $F_0$. Let $L$ be a finite  
extension  of $F_0$. Let $C_L$ be the curve $C\otimes _{K_0} L$ 
obtained by base change and $E_L$ be the corresponding $G_L$-bundle 
over $C_L$.
Then ${\rm Ideg}_{G_L}(E_L)={\rm Ideg}_G(E)$.
\end{lemma}
\proof Let ${\ov F}_0$ be a fixed  algebraic closure of $F_0$. 
Let $E_{{\ov F}_0}$ be the principal $G_{{\ov F}_0}$-bundle over 
$C_{{\ov F}_0}$.  
One first observes that if $(P,\sigma)$ is a pair with the property that 
$P$ is maximal among
all the parabolic subgroup $P'$ containing $B$ 
for which there are reductions $\sigma '$
satisfying 
${\rm Ideg}_{G_{{\ov F}_0}}(E_{{\ov F}_0})={\rm deg}(E_{{\ov F}_0,\sigma '})$ 
then the pair 
 $(P,\sigma)$ defines a Harder-Narasimhan reduction (see Behrend 
\cite{Behrend}). Now by uniqueness of the Harder-Narasimhan 
filtration and the Galois descent argument that
This parabolic reduction is defined over $F_0$. 
If $M$ is a Line bundle over $C$ 
and if $M_L$ is its pull back over $C_L$ then we have 
${\rm deg}(M_L)={\rm deg}(M)$, hence the instability degree does
not change when we take a finite field extension.$\hfill \square$

Let ${\cal E}$ be a principal $G$-bundle over ${\cal C}\ra S$ as above.
We will denote by $x$ a finite field valued point 
$x:{\rm Spec} (k(x))\ra S$, and by $q_x=|k(x)|$. Here $k(x)$ is not 
necessarily the residue field but a finite extension of it. 
Here we fix our notation for $x$. Whenever we say $x$ we mean a finite field 
valued point of $S$.
We also denote by $\E _x$ the principal $G _x$-bundle over the  
curve $\Cu _x$ which is the pullback of the corresponding objects over $S$ 
via $x$. The property we have for $\Cu$ ensures that $\Cu _x$ is a 
smooth projective geometrically connected curve over $k(x)$.

We need the following lemma.

\begin{lemma}\label{bdd}
There exists an integer $N$ such that for each $x$ as above we have 
${\rm Ideg}(\E _x)\le N$.
\end{lemma} 
\proof
Let ${\rm ad}(\E)$ be the adjoint bundle of $\E$. For an $x$ 
$x$  and  a parabolic reduction $(P,\sigma)$ of $\E_x$  we have an 
inclusion ${\rm ad}((\E _x))_{P,\sigma}\hra {\rm ad}(\E _x)$ hence by 
Riemann-Roch Theorem it is enough to bound the degree of 
${\rm ad}(\E _x)$ independent of $x$. But this follows from the 
semi-continuity theorems.$\hfill \square$

Let $G$ be an semisimple simply connected split algebraic group scheme 
over $S$. Let $B$ be a fixed  Borel subgroup and $T$ be a maximal torus 
contained in $B$.
Let ${\cal X}_*(T)$ be the group of one parameter subgroups of $T$. 
Let $[\sigma] \in {\cal X}_*(T)$ be a numerical type. 
Let $E$ be a principal $G$-bundle over $\Cu$. Then the scheme
${\rm Sec}^{[\sigma]}_{\E/B}$ is quasiprojective  over $S$. 
By Corollary \ref{smooth} it follows that ${\rm Sec}^{[\sigma]}_{\E/B}$ 
is smooth over $S$ if  $[\sigma]$ is a numerical type satisfying the property
 $(*)$ for $N>2g-2$. 

We now want to define the Eisenstein series for $\E$. For this we need a 
fixed Borel reduction of $\E$ over all of $S$. By Corollary \ref{ds}
there is a $B$-reduction of $\E$ with some numerical type $[\sigma _0]$
by going to a surjective \'etale extension of $S$. Hence by we may assume that 
our affine scheme $S$ is such that there is a section of the map 
${\rm Sec}^{[\sigma _0]}_{\E/B}\ra S$ for $[\sigma _0] \in {\cal X}_*(T)$.
We fix such a Borel reduction $\sigma _0$ once and for all. 
This has the implication that for any choice of $x$ we have a $B_x$ 
reduction $\sigma _{0,x}$ of 
$\E_x$ of numerical type $[\sigma _0]$.
Now for any $B_x$ reduction of $\E _x$ of numerical type $[\sigma]$ we see 
that $[\sigma]-[\sigma _0] \in {\cal X}_*(T)$  has the property that it lies 
in ${\hat Q}$ (by Lemma \ref{welldefined}). Hence if $\{w_{\alpha}\}$ are the 
fundamental dominant weights then $ ([\sigma]-[\sigma _0])(w_{\alpha})$ 
are integral. Recall from Harder \cite{Harder} the definition of 
$d_{\alpha}([\sigma])= ([\sigma]-[\sigma _0])(w_{\alpha})$.
For simplicity we write $p[\sigma]=[\sigma]-[\sigma _0]$.
Let $\gamma _x(p[\sigma])$ be the cardinality of ${\rm Spec}(k(x))$ 
rational points of $ {\rm Sec}^{p[\sigma]}_{\E _x/B_x}$.
Recall the definition of the Normalized Eisenstein Series 
$$
E(x,\E,\tau)=\sum_{[\sigma] \in {\cal X}_*(T)}\gamma _x(p[\sigma])
q_x^{-\sum_{\alpha \in \Delta}d_{\alpha}([\sigma])}
\Pi_{\alpha \in \Delta}\tau _{\alpha}^{d_{\alpha}([\sigma])}
$$
Here $\tau _{\alpha}$, for $\alpha \in \Delta$ are being thought as variables.
We have suppressed the dependence of the above series on the Borel subgroup 
$B$ and the reduction $\sigma _0$.
Now it follows from Lemma \ref{bdd} that there exists an $N_0$ such that 
for all $x$, $\gamma _x(p[\sigma])=0$ when $d_{\alpha}([\sigma])<-N_0$.
Now we state some of the basic properties of this series which is proved in 
\cite{Harder}

\begin{proposition}\label{rationality}
The Laurent series $E(x,\E,\tau)$ is a rational function on the variables
$\tau$. Moreover it can be written as $E(x,\E,\tau)=P(x,\E,\tau)/Q(x,\tau)$,
where $P$ is a polynomial in $\tau _{\alpha}$'s and $ \tau _{\alpha}^{-1}$'s
for $\alpha \in \Delta$. And 
$$
Q(x,\tau)=\Pi_{\gamma\in \Phi ^+}
(1-q_x\tau _{\gamma})\Pi_{i=1}^{2g}(1-w_i(x)q_x^{-1}\tau _{\gamma}),
$$
where $w_i(x)$ are the eigen values of the Frobenius $Fr _x$ 
on the first cohomology of the curve $\Cu _x$, and 
$\tau _{\gamma}=\Pi_{\alpha \in \Delta}\tau _{\alpha}^{\nu _{\alpha}}$ where
$\gamma =\sum _{\alpha \in \Delta} \nu _{\alpha}w_{\alpha}$.
\end{proposition} 

This proposition  is a consequence  of Theorem 1.6.6  of \cite{Harder}
for  case  ${\underline  x}=1$   and  $\omega  =1$.   From  the  above
proposition  it follows  that  the polynomial  $P(x,\E,\tau)$ has  its
negative powers of  $\tau _{\alpha}$ bounded by $N_0$  for each $\alpha$
as it holds for the Eisenstein series $E(x,\E,\tau)$.

There is a second part of the Theorem 1.6.6 of \cite{Harder} about the 
functional equation satisfied  by the  Eisenstein series using  which we  
get the following upper bound for the degree of the polynomial
$P(x,\E,\tau)$. 
\begin{proposition}\label{ubdd}
There exists a constant $N_1$ such that for each $x$ we have 
$$
-N_0 \le {\rm deg}_{\alpha}P(x,\E,\tau)\le N_1
$$
\end{proposition}
The above proposition is essentially Theorem 1.6.10 of \cite{Harder}. One observes that the constants $N_0$ and $N_1$ depend on the Borel $B$ and the 
reduction $\sigma _0$.
We write the polynomial 
$P(x,\E,\tau)=\sum _{{\underline d}}a(x,{\underline d})\tau^{{\underline d}}$
where ${\underline d}=\{d_{\alpha}\}_{\alpha \in \Delta}$ and 
$\tau^{{\underline d}}=\Pi_{\alpha \in \Delta}\tau _{\alpha}^{d_{\alpha}}$.

The next Lemma we need is about the bound for the coefficients 
$a(x,{\underline d})$. 

\begin{lemma}\label{points}
There exists a constant $M$ independent of $x$ such that for each $x$ we have 
$|a(x,{\underline d})|\le q_x^M$.
\end{lemma}
\proof
For a fixed $[\sigma]$, we know that $ {\rm Sec}^{p[\sigma]}_{\E/B}$ is a finite type quasi-projective scheme over $S$ hence  
we can find an constant $M_{p[\sigma]}$ such that for each $x$ we have a bound
$\gamma _x(p[\sigma])\le q_x^{M_{p[\sigma]}}$. 
It follows from Lemma \ref{lowerbound} that there are only finitely many 
$[\sigma]$ with the property that 
\begin{equation}\label{2bound}
-N\le d_{\alpha}([\sigma])\le N_1 
\end{equation}
for each $\alpha$. Hence we can find a single constant $M_0$ such that 
$\gamma _x(p[\sigma])\le q_x^{M_0}$ for all $[\sigma]$ satisfying 
(\ref{2bound}).
Now we use the Proposition \ref{rationality} to write
$E(x,\E,\tau)=P(x,\E,\tau)/Q(x,\tau)$. Expanding both sides of the series 
 using the power series expansion of $1/Q(x,\tau)$, 
the coefficients $a(x,{\underline d})$ can be computed as a 
linear combinations of $\gamma _x(p[\sigma])$ for $[\sigma]$ satisfying 
(\ref{2bound}) (by inverting an upper triangular matrix) 
with coefficients as powers of $q_x$ and $w_i(x)$. This proves
the Lemma.$\hfill \square$

Now we follow the arguments of the Theorem 2.3.1 and Theorem 2.3.2 of 
\cite{Harder} we get the following.
\begin{proposition}\label{hard}
There exists a constants  $N$ and $C$ 
independent of $x$ such that if $[\sigma]$ satisfies the property $(*)$ 
for $N$ then we have 
$$
|\gamma _x([\sigma])|= q_x^{d([\sigma])+(1-g){\rm dim}(G/B)}+
R_x([\sigma])
$$ 
where 
$$
R_x([\sigma]) \le Cq_x^{d([\sigma])+(1-g){\rm dim}(G/B)-1/2}
$$ 
\end{proposition}
The basic idea of proof of the above proposition is 
that we know explicitly the poles of the Eisenstein series namely the 
zeros of the polynomial $Q(x,\tau)$. The next step is the 
explicit computation of the residue of the Eisenstein series at the point 
$\{\tau _{\alpha}\}=\{ 1/q_x\}$. Now one writes the series  
$E=E_1(1/\Pi _{\alpha \in \Delta}(1-q_x\tau_{\alpha}))+E_2$
where $E_1$ is essentially the residue of $E$. Then one observes that the series $E_2$ has better radius of convergence.
From here it follow that for large $N$ the coefficients of the Eisenstein 
series are dominated by the coefficients of $E_1$ and the explicit residue computation now yield the Proposition.

Now by applying the results of Lang and Weil \cite{Lang-Weil} we see that 
if $[\sigma]$ satisfies $(*)$ for $N$ then after a finite base change of $x$, 
${\rm Sec}^{[\sigma]}_{\E _x/B}$ has a unique 
irreducible component of maximal dimension and this dimension is equal to
$d([\sigma])+(1-g){\rm dim}(G/B)$. 
By deformation theoretic lower bounds we see that every component is 
atleast of this dimension. Hence we see that 
${\rm Sec}^{[\sigma]}_{\E _x/B}$ is absolutely irreducible for each $x$.

Now we have a finite type smooth morphism ${\rm Sec}^{[\sigma]}_{\E /B}\ra S$. 
By taking an affine open subscheme $S_0={\rm Spec}(A)$ of $S$ we may assume 
that there is 
a dense affine open subscheme $U={\rm Spec}(R)$ of 
${\rm Sec}^{[\sigma]}_{\E /B}\ra S$ 
such that $U$ surjects onto $S_0$. Hence we obtain a finite type 
faithfully flat morphism $U\ra S_0$ with the property that for each $x$ in 
$S_0$  the fiber $R\otimes _A k(x)$ is irreducible.
To show that ${\rm Sec}^{[\sigma]}_{\E _K/B}$ is irreducible it is enough to show that $R\otimes _A K$ has no non trivial (not equal to $0$ or $1$) 
idempotent elements  as $R$ is a smooth $A$ algebra. If $R\otimes _A K$ has 
an idempotent element then there is an 
$f\in A$ such that $R_f=R\otimes _A A_f$ contains a non-trivial idempotent 
element.  Hence we may assume by replacing $A$ by $A_f$ that $R$ 
contains a non trivial idempotent element $e$. 
Since $R$ is reduced, there is an open subscheme $U_1$ of 
${\rm Spec}(R)$ such that for 
${\mathfrak p}\in U_0$ the image of $1-e$ in 
$R_{\mathfrak p}/{\mathfrak p}R_{\mathfrak p}$ is non zero. Since $R$ is 
smooth over $A$, the image of $U_1$  
in ${\rm Spec}(A)$ is open. Hence by shrinking $A$ to another $A_f$ we may 
assume that the non-trivial idempotent 
element $e\in R$ has the property that the image of $1-e$ in 
$R\otimes _A k(x)$ is non-zero for each $x$. 
Now by irreducibility of 
$R\otimes _A k(x)$ we see that $e$ maps to $0$ for each $x$. Hence  
$e\in \bigcap {\mathfrak m}R$ where the intersection is over all maximal ideals of $A$.

\begin{lemma}Let $A$ be a finite type algebra over $F$ 
($=\Z$ or $\F_q$). If $R$ is finite type over $A$ then 
$\bigcap {\mathfrak m}R ={\rm rad}(R)$.
\end{lemma}
\proof
Firstly the conditions of the lemma ensures that $R$ is a Jacobson ring hence 
${\rm rad}(R)$ is the intersection of all maximal ideals of $R$.  Again since $R$ is finitely generated over $F$ ensures that for every maximal ideal 
${\mathfrak n}$ in $R$ the field $R/{\mathfrak n}$ is finite. Hence we see that the 
morphism ${\rm Spec}(R) \ra {\rm Spec}(A)$ takes maximal ideals to maximal ideals.
of $A$. This implies the lemma. $\hfill \square$

Now by above lemma, since $R$ is reduced, it follows that $e=0$. 
Hence the proof of the theorem \ref{borel} is complete. $\hfill \square$

\rem Note that we need to actually show that 
${\rm Sec}^{[\sigma]}_{\E _K/B}$ is absolutely irreducible. For this we have to further show that if $L$ is a finite separable extension of $K$ then
${\rm Sec}^{[\sigma]}_{\E _K/B}$ remains irreducible. We take the 
normalization $A'$ of $A$ in $L$. Then the whole setup pulls back to the 
setup over $A'$. Since ${\rm Spec}(A') \ra {\rm Spec}(A)$ takes maximal ideals 
to maximal ideals hence above proof actually shows that for  
the same $N$, the Proposition \ref{hard} works in the setup for $A'$.

We return to our notations where the objects are over an algebraically closed field $k$. Now we extend the above results to the case of parabolic subgroups.

We need the following lemma about the behavior of the Harder-Narasimhan 
reduction of a principal $G$-bundle which vary over a family.

\begin{lemma}\label{canred}
Let ${\cal E}\ra C\times S$ be a family of $G$-bundles with $S$ an
irreducible scheme of finite type over $k$. Then there is a non-empty
open subset $U$ of $S$, a parabolic $P$ and reduction $\sigma _u$ of
${\cal E}|_{C\times \{u\}}$ to $P$ for each $u\in U$ such that 
$\sigma _u$ is the Harder-Narasimhan reduction of ${\cal E}|_{C\times \{u\}}$ for each
$u\in U$ and the numerical types $[\sigma _u]=[\sigma _v]$ for every
$u,\,v\,\in U$.
\end{lemma}
\proof
Consider the family ${\rm ad}({\cal E})$ over $C\times S$, 
since the degrees of the subbundles of each of the vector bundles 
occurring in the above family is bounded above by an integer
which depends only on $S$, there are only finitely many choices 
for the $(P,[\sigma])$ with $P\supset B$ where the $G$-bundles in the 
family ${\cal E}$
can admit a Harder-Narasimhan reduction to these parabolics with 
numerical types $[\sigma]$.
Using the uniqueness of the Harder-Narasimhan reduction we  
produce a finite collection of constructible subsets of $S$ (by Lemma 
\ref{lclosed})
corresponding to $G$-bundles in the above family whose Harder-Narasimhan
 reduction is defined by the pair $(P,[\sigma])$. 
The union of these constructible 
subsets is $S$. Hence one of them contains a non-empty open subset of $S$. 
This proves the lemma.  $\hfill \square$

Let $P$ be a fixed parabolic subgroup of $G$ containing $B$. 
Let $L$ be the Levi quotient of $P$. We denote by ${\ov L}$ the quotient 
$L/Z^0(L)$ where $Z^0(L)$ is the connected component of the center of $L$.
Let ${\ov B}$ (resp. ${\ov T}$) be the Borel subgroup (resp. maximal torus)
be defined by image of $B$ 
(resp. $T$). This defines an induced root system for ${\ov L}$.
Let ${\hat Q}_L$ be the coroot lattice of this root system.
Let $\sigma$ be a $P$ reduction of $E$ with the numerical type 
$[\sigma]$. We will denote by
$c_L([\sigma])$ the element $c_L(p_*(E_{\sigma}))$. This is well defined by 
Lemma \ref{G-L}.
Recall the definition of $M_G(c,d)$ from the paragraph above Proposition 
\ref{bound}.

\begin{proposition}\label{important}
Let $E$ be a principal $G$-bundle over $C$.
There exists a constant $D$ with the property if $E$ admits a $P$ reduction 
of numerical type $[\sigma]$ then the  
scheme of sections ${\rm Sec}^{[\sigma]}_{E/P}$ has an open dense subscheme 
$U^{[\sigma]}$ with the property that the ${\ov L}$-bundles associated to 
every point of $U^{[\sigma]}$ is a member of $M_{{\ov L}}(c_L([\sigma]),D)$. 
\end{proposition}
\proof 
Let $X$ be an irreducible component of ${\rm Sec}^{[\sigma]}_{E/P}$.
Let $\E _P$ be the restriction of the  universal $P$-bundle over $C\times X$
Let $p$ be the natural surjection from $P$ to ${\ov L}$. 
Then we have a family $p_*\E _P$ of ${\ov L}$-bundles over $C \times X$. 
By Lemma \ref{canred} there is a parabolic subgroup 
${\ov P_1}$ of ${\ov L}$ and an open subscheme $U$ of $X$ such that  
$(p_*\E _P)_x$ admits Harder-Narasimhan reduction ${\ov \sigma}_{1,x}$ to 
the parabolic ${\ov P_1}$ with a fixed  numerical type 
$[{\ov \sigma}'_x]=[{\ov \sigma}_1]$ for each $x\in U$.
Let $P_1$ be the parabolic subgroup of $G$ contained in $P$ 
whose image in ${\ov L}$ is ${\ov P_1}$. 
Now by Lemma \ref{compare} the reduction of structure group 
${\ov \sigma}_{1,x}$ 
canonically defines a reduction of structure group $\sigma _{1,x}$ of 
$(E_P)_x$ to $P$ with 
the property that $[{\sigma}_{1,x}]=[\sigma 1]$ for each 
$x\in U$ and such that  $[\sigma _1]$ maps to 
$[\sigma]$ (resp.  $[{\ov \sigma _1}]$) under the homomorphism
${\cal X}_*(P_1)\ra {\cal X}_*(P)$ 
(resp. ${\cal X}_*(P_1)\ra {\cal X}_*({\ov P_1})$ is $[{\ov \sigma _1}]$).
By definitions the instability degree ${\rm Ideg}(p_*(\E_P)_x)$ at $x\in U$
 is exactly computed to be equal to $-d([{\ov \sigma}_1])$. Using the exact 
sequence 
$$ 
0\ra T_{{\ov \sigma}_{1,x}}\ra T_{\sigma _{1,x}}\ra T_{\sigma}\ra 0  
$$
we conclude that $d([{\ov \sigma}_1])= d([\sigma _1])-d([\sigma])$.

Now we have a natural morphism 
${\rm Sec}^{[\sigma _1]}_{E/P_1} \ra {\rm Sec}^{[\sigma]}_{E/P}$
whose image contains $U$.
This implies that ${\rm dim}(X)\le {\rm dim}({\rm Sec}^{[\sigma _1]}_{E/P_1})$.
Now we use the deformation theoretic lower bounds for the dimension of $X$ and the upper bound Theorem \ref{stable} we obtain for each $x \in U$
$$
{\rm Ideg}(p_*(\E_P)_x) \le (g-1){\rm dim}(G/P) + {\rm dim}(G/P_1) 
- {\rm Min}_{i=1\cdots m}\{d(\gamma _i)\}
$$
Where $\gamma _i$'s are the minimal numerical types for $E$ with respect 
to $P_1$.
Now the proposition follows as the right hand side of the above expression is 
independent of $[\sigma]$. $\hfill \square$

We will fix the constant $D$ prescribed by the above proposition 
once and for all. 
This will allow us to work with the open dense subscheme $U^{[\sigma]}$ 
of  ${\rm Sec}^{[\sigma]}_{E/P}$.

We also need the following lemma
Let $c\in {\cal X}_*(T)/{\hat Q}$ be fixed. 
\begin{lemma}\label{reductions}
For every positive integer $n$ there exists an integer $M_c(n)$ such that 
every member in $M_G(c,D)$ admits a reduction to $B$ with numerical type 
$[\sigma]$ satisfying $n\le [\sigma](\alpha)\le M_c(n)$ for each 
$\alpha \in \Delta$.
\end{lemma}

\proof This a an immediate consequence of the corollary \ref{ds} where we 
just choose $M_c(n)$ to be  the maximum of $[\sigma](\alpha)$ for 
$\alpha \in \Delta$. $\hfill \square$

The main Theorem of the section is the following.
\begin{theorem}\label{main5}
Let $E$ be a principal $G$-bundle over $C$. Let $P$ be a 
parabolic subgroup of $G$. 
There exists an integer $N$ such that 
if $E$ admits a $P$ reduction of numerical type $[\sigma]$ satisfying the 
property $(*)$ for $N$ then 
${\rm Sec}^{[\sigma]}_{E/P}$ is an irreducible and generically smooth 
of dimension $d([\sigma])+(1-g){\rm dim}(G/P)$.
\end{theorem} 
\proof
Let $I$ be the subset of $\Delta$ which defines the parabolic subgroup $P$.
Hence the simple roots of the quotient ${\ov L}$ are exactly the simple roots 
in $I$.
If $\beta \nin I$ is a  character of $T$ then $\beta  |_{Z^0(L)}$ is 
nontrivial. Hence there is an positive integer $n_{\beta}$ such that 
 $n_{\beta}\beta$ extends to a character of $L$ (hence $P$). This implies that 
there exists a character $\chi _{\beta}$ of $P$ such that 
\begin{equation} \label{eqa}
\chi _{\beta}|_T=n_{\beta}\beta +\sum_{\alpha \in I}n_{\beta,\alpha}\alpha.
\end{equation}
We note that if $\chi$ is any non-trivial character of $L$ which when 
restricted to $T$ is a non-negative linear combination of simple roots 
then some positive integral multiple of $\chi$ is a non-negative linear 
combination of $\chi _{\beta}$ for $\beta \nin I$.
Let $m_I= {\rm Max}_{\beta \nin I}\{\sum _{\alpha \in I}|n_{\beta,\alpha}|\}$
Hence if $[\sigma]$ satisfies $(*)$ for $N$ then 
$[\sigma](\chi _{\beta})\ge N$ and conversely if 
$[\sigma](\chi _{\beta})\ge N$ for each $\beta \nin I$ then 
$[\sigma]$ satisfies $(*)$ for $N_1$ where $N_1= N/n_I$ and 
$n_I={\rm Max}_{\beta \nin I} \{n_{\beta}\}$.

By Theorem \ref{borel} there exists an integer $N_B$ such that 
for each $[\sigma _0] \in  {\cal X}_*(T)$ satisfying 
$(*)$ for $N_B$ the space of sections ${\rm Sec}^{[\sigma]}_{E/B}$
is irreducible and smooth of dimension $d([\sigma _0])+(1-g){\rm dim}(G/B)$.

Let $N_P$ be defined by $N_P=n_I N_B+m_IM_D$, where $D$ is prescribed 
by the Proposition \ref{important} and 
$M_D ={\rm Max}_{c\in {\cal X}_*({\ov T})/Q_L}\{M_c(N_B)\}$. 
Here $ M_c(N_B)$'s are the constants prescribed by the Lemma  
\ref{reductions} for the choice of $D$ and $n=N_B$.
The constant $M_D$ is finite because ${\ov L}$ is semisimple 
algebraic group, hence  there are only finitely many classes in
${\cal X}_*({\ov T})/Q_L$.

Let $[\sigma] \in {\cal X}_*(P)$ be such that $[\sigma]$ satisfies $(*)$ for 
$N_P$ then we claim that  ${\rm Sec}^{[\sigma]}_{E/P}$ satisfies the properties mentioned in the statement of the Theorem.
Infact if $U^{[\sigma]}$ is the open dense subscheme prescribed by the 
Proposition \ref{important} then we show that $U^{[\sigma]}$ is 
irreducible and smooth of expected dimension.

Let $\E$ be the family of $P$-bundles over $C\times U^{[\sigma]}$ defined by the universal properties of the space of sections.
Then by Lemma \ref{reductions}
for each $x\in U^{[\sigma]}$ the ${\ov L}$-bundle $p _*\E_x$ admits 
${\ov B}$ reductions with numerical type $[{\ov \sigma}_1]$ satisfying
$N_B \le [{\ov \sigma}_1](\alpha) \le M_D$.
Now by Lemma \ref{compare} these reductions canonically defines $B$-reductions 
of $\E_x$ with numerical type $[\sigma _1]$ such that  
the image of $[\sigma _1]$ under the map ${\cal X}_*(B)\ra {\cal X}_*(P)$
(resp. ${\cal X}_*(B)\ra {\cal X}_*({\ov B})$) is $[\sigma]$ 
(resp. $[{\ov \sigma} _1]$).
Hence $U^{[\sigma]}$ lies in the image of the natural morphism 
$f:{\rm Sec}^{[\sigma _1]}_{E/B} \ra {\rm Sec}^{[\sigma]}_{E/P}$.

The Theorem now follows if we show that $[\sigma _1]$ satisfies the property
$(*)$ for $N_B$. This is because the irreducibility of 
${\rm Sec}^{[\sigma _1]}_{E/B}$ would imply irreducibility of $U^{[\sigma]}$.
Moreover for a reduction $\sigma \in U^{[\sigma]}$ if 
$\sigma _1 \in {\rm Sec}^{[\sigma _1]}_{E/B}$  is such that 
$f(\sigma _1)=\sigma$ then the $B$-bundle $E_{\sigma _1}$ extends to 
$E_{\sigma}$. Since the map $H^1(C,\,T_{\sigma _1}) \ra 
H^1(C,\,T_{\sigma})$ is surjective, the fact that $H^1(C,\,T_{\sigma _1})=0$
implies that $U^{[\sigma]}$ is smooth at $\sigma$.

Hence we have to only show that $[\sigma _1](\alpha)\ge N_B$ for each 
$\alpha \in \Delta$.
If $\alpha \in I$ then already we have 
$[\sigma _1](\alpha) =[{\ov \sigma}_1](\alpha) \ge N_B$.
If $\beta \nin I$ then the inequality 
$[\sigma _1](\beta) \ge N_B$ follows from the 
above inequality for $\alpha \in I$ and the Equation (\ref{eqa}) 
by observing that  
$[\sigma _1](\chi _{\beta}|_T)=[\sigma ](\chi _{\beta})\ge N_P$. 
This completes the proof of Theorem \ref{main5}. $\hfill \square$
\rem The results of this section holds for any connected reductive 
algebraic groups. To do this firstly we may assume that the group 
$G$ has no center as the
map $G\ra G/Z$ has the property that $P$-reduction of a principal
$G$-bundle $E$ is in bijection with $P/Z$-reductions of the
associated $G/Z$-bundle. This reduces the problem to semisimple algebraic 
group. Let $f:{\wt G} \ra G$ be the simply connected cover. 
Then one checks that all the proofs go through if
we take semisimple simply connected group schemes over the curve instead of 
the principal bundles.  Now one proves that for any principal $G$-bundle
there is a ${\wt G}$ group scheme whose quotient modulo center in the flat 
topology is the group scheme over $C$ associated to the $G$-bundle $E$.
The last statement follows from the existence of a connected reductive 
algebraic group $H$, and a surjective homomorphism (in the fppf topology)
$g:H\ra G$ such that ${\rm ker}(g) =(\G_m)^n$ and $[H,\, H]={\wt G}$.
This is because  $E$ comes from a 
$H$-bundle $E_1$ (as $H^2(C,\, \G_m)=0$) hence 
the group scheme we are interested is the associated 
fiber space $E_1({\wt G})$ for the conjugation action of $H$ on ${\wt G}$. 
Since there exists an embedding of the group scheme 
${\rm ker}f \subset (\G_m)^n$ for some $n$, we check that the  quotient 
$({\wt G}\times (\G_m)^n)/{\rm ker}f$ (in the fppf topology) for the 
diagonal action of ${\rm ker}f$ defines such a choice of $H$.

\section{Generic stable bundles}

In this section we prove the main results about the structure 
of the moduli spaces of stable maps ${\ov M}_g(E/P,\,\beta _{[\sigma]})$ under 
stronger assumptions on the principal $G$-bundle $E$.
In this section we will also assume that the genus of the curve $C$ is 
atleast two. This will ensure that there are stable $G$ bundles of every 
topological type by Proposition \ref{exis}. 

Following the Example 5.7 of \cite{Popa-Roth} we define the notion of a 
generically stable $G$-bundles as follows.

\begin{definition}
We say a principal $G$-bundle $E$ is generically stable if for any parabolic $P$ and a $P$ reduction $\sigma$ of $E$ the following holds.
\begin{equation}\label{def}
{\rm dim}({\ov M}_g(E/P,\,\beta _{[\sigma]}))=d([\sigma])+(1-g){\rm dim}(G/P) 
\end{equation}
\end{definition}
Note from the definition that generically stable is actually stable.
Our main result is the following.

\begin{theorem}\label{generic}
Let $C$ be a curve of genus atleast two. Let $c \in {\cal X}_*(T)/{\hat Q}$
be fixed. Then there exists a generically stable $G$-bundle with topological 
type $c$.
\end{theorem}
\proof 
Let $c \in {\cal X}_*(T)/{\hat Q}$ be fixed. Let $\E$ be a family of 
$G$-bundles on the curve $C$ parameterized by a finite type smooth scheme $S$  
which is miniversal at every point of $S$ and such that 
$\E _x=\E |_{C\times \{x\}}$ is stable for each $x \in S$. Such a family 
exists by Proposition \ref{versal} and \ref{exis}.

Let $\Gamma_1$ be the set of numerical types satisfying the property that
$d([\sigma])\le g\cdot {\rm dim}G/P$ and there exists a 
$y\in S$ such that $\E _y$ admits a $P$ reduction with numerical type
$[\sigma]$. By lemma \ref{lowerbound} and  Lemma \ref{bdd} it follows that 
$\Gamma_1$ is a finite set.

\begin{lemma}\label{openx}

There exists an non-trivial open subset 
$U_P \subset S$ with the property that if for some $x\in U_P$ the $G$-bundle 
$\E _x$ admits a $P$ reduction of numerical type $[\sigma]\in \Gamma $ 
then for every $y \in U_P$ the bundle $\E _y$ admits a $P$ reduction of 
numerical type $[\sigma]$.
\end{lemma}
\proof

Let $[\sigma] \in \Gamma _1$ be an numerical type. 
Let $V_{[\sigma]}$ be the constructible subset of 
$S$ consisting of bundles in $\E _y$  which admit $P$ reductions of 
numerical type $[\sigma]$ (by Lemma \ref{lclosed}).
Let $\Gamma \subset \Gamma_1$ be the subset of numerical types 
with the property that 
for each $[\sigma] \in \Gamma$ the constructible set $V_{[\sigma]}$ contains a
non-empty open subset of $S$.

By Theorem 1.1 of Holla-Narasimhan \cite{Holla-Narasimhan},
the union of the 
these finitely many constructible sets $V_{[\sigma]}$ for 
$[\sigma] \in \Gamma _1$
is all of $S$. Hence one of them
must contain a non-empty open subset of $S$. This proves that
$\Gamma$ is non-empty. Let 
$$
U= \bigcap _{[\sigma]\in \Gamma}V_{[\sigma]}\bigcap _{[\sigma _1]\in 
\Gamma_1-\Gamma} ({\ov V}_{[\sigma _1]})^c
$$ 
where $({\ov V}_{[\sigma _1]})^c$ is 
the complement of the closure of $V_{[\sigma _1]}$.
Then $U$ contains a non-empty open subset of $S$ which satisfies the properties mentioned in the lemma. $\hfill \square$

Let $P$ be a parabolic subgroup containing $B$. Let $\Gamma$ be the finite 
set numerical types for the family $\E$ as described before.  
By Lemma \ref{openx} we have an open subscheme $U_P$ of $S$ satisfying the properties mentioned in the Lemma.
Now define $U=\bigcap U_P$ where the intersection is over the 
finitely many parabolic subgroups containing the Borel subgroup $B$.
Hence $U$ is a non-empty open subscheme of $S$ with the property that 
the defining morphism $f:{\rm Sec} ^{[\sigma]}_{{\cal E}|_{U}/P} \ra U$ is 
surjective for each $P\supset B$.
We will show that for any $x \in U$ the principal bundle $E=\E _x$ is 
generically stable.

Let $\sigma$ be a $P$ reduction of $E$ such that $[\sigma]\in \Gamma$.
To simply notations we will denote by $Y$ the space 
 ${\rm Sec} ^{[\sigma]}_{{\cal E}|_{U}/P}$. 

If $y\in Y$ is the point corresponding to $\sigma$ then  
${\rm Sec} ^{[\sigma]}_{E/P}$ is the fiber $Y_{x}$ of the morphism 
$f:Y\ra U$ at the point $x=f(y)\in U$.
We denote by $df:T_yY\ra T_{x}U=H^1(C,\,{\rm ad}E)$ 
the induced map on the tangent space (Here the schemes are not smooth and by tangent space we mean the vector space of $k[\epsilon]/\epsilon ^2$ valued 
points whose associated $k$-valued point is $y$). Then
we have the following exact sequence
\begin{equation}\label{exact}
T_y{\rm Sec} ^{[\sigma]}_{E/P}= 
H^0(C,\, T_{\sigma})\ra T_yY \ra H^1(C,\,{\rm ad}E), 
\end{equation}
where $T_{\sigma}$ is the pull back of the tangent bundle of the fibers of the 
morphism $E/P \ra C$ by $\sigma$.
Let ${\cal V}\lra C\times S_P$ be a 
family of $P$-bundles, parameterized by finite type smooth scheme $S_P$,
which is miniversal at the point $y'=E_{\sigma}\in S_P$.
Now the family of $P$-bundles
defined by $S_P$ when extended to $G$ defines a family of $G$-bundles. 
There is a stable $G$-bundle in this family namely $E$ which is an 
extension of $E_{\sigma}$. Hence there is a non-trivial open subset 
$S_0$ of $S_P$ which
corresponds to points which extend to stable $G$-bundles (by Proposition 
\ref{open}). Hence by replacing $S_P$ by an \'etale neighborhood  
we may assume there is a morphism 
$S_P \lra U$.
Consider the family ${\cal E}_P$ of $P$-bundles on $C\times Y$ obtained by 
the universal property of $Y$.
Now by going to an \'etale neighborhood $V$ of $y \in Y$  and of $S_P$,
and an automorphism of $U$ in a neighborhood of $x \in U$
we may assume that the restriction of $f$ (again denoted by $f$) defines a 
morphism $f:V \ra U$ which can be written as $j \circ g$, where 
$g:V \ra S_P$ is defined by the versal property of $S_P$ at $E_{\sigma}$ 
and $j:S_P \ra U$ by the versal property of $U\subset S$ at $E$.

Now this implies that the map $df:T_yV\ra H^1(C,\,{\rm ad}E)$ factors
through $T_{y'}S_P=H^1(C,\,{\rm ad}E_{\sigma})$. 
Hence we have the following commuting diagram which has exact horizontal rows
$$
\begin{array}{ccccccc}
0                           & \lra                    &
H^0(C,\,T_{\sigma})         & \lra                    & 
T_yV                        &  \stackrel{df}{\lra}    &
H^1(C,\,{\rm ad}E)     \\ 
\downarrow                  &                        &
\parallel                   &                        &
\downarrow dg               &                        &
\parallel              \\              
H^0(C,\,{\rm ad}E)          & \stackrel{\eta}{\lra}   & 
H^0(C,\,T_{\sigma})         & \stackrel{\delta}{\lra} &  
H^1(C,\,{\rm ad}E_{\sigma}) & \stackrel{dj}{\lra}     & 
H^1(C,\,{\rm ad}E) 
\end{array}
$$

From the above diagram it follows that 
${\rm ker}(dg) = {\rm im}(\eta)$. This gives us a dimension bound
${\rm dim}(T_yV)\le {\rm dim}(H^1(C,\,{\rm ad}E_{\sigma}))
+{\rm dim}({\rm im}(\eta))$. Using the following exact sequence 
$$
0 \lra H^0(C,\,{\rm ad}E_{\sigma})\lra 
H^0(C,\,{\rm ad}E)\stackrel{\eta}{\lra}  H^0(C,\,T_{\sigma}) 
$$
and Riemann-Roch Theorem for $E_{\sigma}$ we get the dimension bound 
${\rm dim}(T_yV)\le d([\sigma])+(g-1){\rm dim}(P)
+{\rm dim}(H^0(C,\,{\rm ad}(E)))$. 
The Equation (\ref{exact}) now implies that
\begin{equation}\label{Xl}
{\rm dim}({\rm Sec} ^{[\sigma]}_{E/P})\le d([\sigma])+(1-g){\rm dim}(G/P).
\end{equation}

But this is exactly the deformation theoretic lower bound estimate of
the dimension of the space of sections. Hence we have proved that the 
inequality
(\ref{Xl}) is an equality.  This proves the Equation (\ref{def}) for 
$[\sigma] \in \Gamma$ assuming there are no pathological components in
${\ov M}_g(E/P,\,\beta _{[\sigma]})$.

For an arbitrary $[\sigma]$. Let $X$ be any irreducible component of 
${\ov M}_g(E/P,\,\beta _{[\sigma]})$. By Remark \ref{veryuseful}, 
there exists a $P$ reduction $\sigma _0$ of $E$ with 
$[\sigma _0] \in \Gamma$ and 
$[\sigma]\le [\sigma _0]$ 
such that 
${\rm dim}(X)\le {\rm dim}(X_0)+d([\sigma])-d([\sigma _0])-k$. 
Here $k$ is the number of nodes in the curve corresponding to a general 
point in $X$ and $X_0$ is an irreducible component of 
${\ov M}_g(E/P,\,\beta _{[\sigma  _0]})$ containing $\sigma _0$.
Now the Theorem  \ref{generic} follows from (\ref{Xl}) for $\sigma _0$ 
and the deformation theoretic lower bounds for the dimension of $X$ 
(also we get $k=0$). $\hfill \square$

\rem Note that even for a stable bundle 
$E$, in general $H^0(C,\,{\rm ad}(E))$ is not equal to the center of the 
Lie algebra of $G$ (this is true in characteristic 0). So, one has to be 
careful at this stage. Also the use of versal family is again due to the 
lack of the moduli space of stable bundles in positive characteristic. Even 
in characteristic zero we need to use the versal family because of 
the non-representability of the functor $H^1_S(C,\,R_uP(\E _L))$ for 
a family of $L$ bundles $\E _L$.

\rem The proof of the above theorem actually show that there exists an 
open subset of stable bundles in any family which are generically stable.
 
\rem  Let $N\ge g\cdot{\rm dim}(G/P)$ be an integer. 
Let $c \in {\cal X}_*(T)/{\hat Q}$ be fixed.
Let $\Gamma$ be the finite set of  
numerical types for $P$ defined by the property that for each 
$[\sigma] \in \Gamma$ there exists a stable bundle to topological type $c$ 
which admits a $P$-reduction 
of numerical type $[\sigma]$ and $d([\sigma])\le N$. 
The proof of the above theorem actually shows that if the genus of 
the curve is atleast two then there exists a stable bundle $E$ 
(hence an open set of stable bundles) such that if $E$ admits a $P$ reduction 
of numerical type $[\sigma] \in \Gamma$ then ${\rm Sec}^{[\sigma]}_{E/P}$ is 
smooth. This is because the in the proof we actually get the dimension bound 
for the tangent space of ${\rm Sec}^{[\sigma]}_{E/P}$. 
Hence the smoothness follows from the deformation theoretic lower bounds.

As a consequence of the Theorem \ref{generic} we get a result which 
generalizes the lower bound theorem of  Lange (Satz 2.2,\cite{Lange}).

\begin{corollary}\label{main3}Let $E$ be a generic stable bundle 
 if $E$ admits a $P$ reduction $\sigma$ then 
we have $d([\sigma])\ge (g-1){\rm dim}(G/P)$.
\end{corollary}
\proof 
This follows from the definition of a generic stable bundle. The main point is 
the existence of such bundles which is the content of the Theorem 
\ref{generic}  

Now we prove a result on the structure of the space of parabolic reductions 
of a generically stable bundle.

\begin{proposition} \label{lgeneric}
Let $E$ be a generic stable $G$-bundle. Let $[\sigma]$ be an numerical type.
Assume that  
${\ov M}_g(E/P,\,\beta _{[\sigma]})$ is non-empty.  Then a
generic element in every component of 
${\ov M}_g(E/P,\,\beta _{[\sigma]})$ 
corresponds to reduction of structure group of $E$ to $P$ with the property 
that the associated Levi $L$-bundle is generically stable.
\end{proposition}
\proof
If $X$ is an irreducible component of 
${\ov M}_g(E/P,\,\beta _{[\sigma]})$ and if the generic element of $X$ 
corresponds to a map from a curve with $k$ nodes then
using the Remark \ref{veryuseful} we see that 
${\rm dim}(X)\le  {\rm dim}(X_0) +d([\sigma])-d([\sigma _0])-k$ for some 
$[\sigma _0]$ and $X_0$ an irreducible component of 
${\ov M}_g(E/P,\,\beta _{[\sigma _0]})$. Since $X_0$ has the expected 
dimension we get a contradiction to the deformation theoretic lower bound
for $X$. This proves that there are no pathological components in 
${\ov M}_g(E/P,\,\beta _{[\sigma]})$.
Hence it is enough to prove the proposition for an irreducible component 
$S$ of ${\rm Sec}^{[\sigma]}_{E/P}$. By the universal property we have a 
family $\E_P$ of $P$ bundles over $C\times S$. Let $\E _{{\ov L}}$ be the 
associated family of ${\ov L}$ bundles over $C\times S$, where ${\ov L}$ is 
the quotient of Levi subgroup of $P$ by its connected component of the center.
For any parabolic ${\ov P}_1$ of ${\ov L}$ there is a parabolic 
$P_1 \subset P$ such that ${\ov P}_1$ is the image of $P_1$ under the map 
$p:P\ra {\ov L}$.

Now we continue as in the proof of the Theorem \ref{generic} with the 
family over ${\ov L}$ instead of $G$. Using the Equation (\ref{def}) one 
checks that the same proof works once we use the following Lemma.

\begin{lemma} With above notations, 
let $[\sigma _1] \in {\cal X}_*(P_1)$ be a numerical type which maps to 
$[\sigma]$  and to a numerical type $[{\ov \sigma} _1]$ of ${\ov P}_1$.  
Then the spaces  
${\rm Sec}^{[\sigma _1]}_{E/P_1}$ and 
${\rm Sec}^{[{\ov \sigma} _1]}_{\E_{\ov L}/{\ov P}_1}$
are naturally isomorphic.
\end{lemma}
\proof  This is a simple consequence of the universal properties of the 
space of sections.  $\hfill \square$ 

We also get the following result which generalizes the Proposition 
\ref{main2} for the case of parabolic subgroups when $E$ is generically stable.
\begin{corollary}
Suppose $[\sigma]$ and $[\tau]$ are two numerical types with the property
$[\tau] \le [\sigma]$. Suppose $E$ is a generic stable $G$-bundle
which admits a reduction of structure group $\sigma$ to $P$ with
numerical type $[\sigma]$. Then $E$ admits a reduction of structure group 
to $P$ with numerical type $[\tau]$.
\end{corollary}
\proof
This immediately follows from the Proposition \ref{lgeneric} once we show 
that the space ${\ov M}_g(E/P,\,\beta _{[\tau]})$ is non-empty and this is so 
because we can always attach rational tails to a map from a curve 
corresponding to reduction $\sigma$.  $\hfill \square$

Now we prove result which generalizes the Theorem 6.7 of 
\cite{Popa-Roth} and gives a characterization of generically stable bundles.
\begin{proposition} \label{main1}
A principal $G$ bundle $E$ is generically stable if and only if there exists 
an integer $N$ such that ${\ov M}_g(E/P,\,\beta _{[\sigma]})$ is irreducible 
for all parabolic subgroups and numerical types $[\sigma]$
satisfying the property {\rm ($*$)} for $N$
\end{proposition}
\proof 
If a $G$-bundle $E$ is generically stable then by Proposition \ref{lgeneric}
it follows that ${\ov M}_g(E/P,\,\beta _{[\sigma]})$ has no pathological 
components. Now the ``only if part''  follows from Theorem \ref{main5}.
For the other way, using the Theorem \ref{main5} again we find an 
$N_1 \ge N$ such that if $E$ admits a $P$ reduction of numerical type 
$[\sigma]$ satisfying ($*$) for $N_1$ then 
${\rm {\ov M}_g(E/P,\,\beta _{[\sigma]})}$ has expected dimension. Hence by 
Remark \ref{useful} we conclude that there is a $P$ reduction $\sigma _0$ of
 $E$ and a component $X_0$ of ${\rm {\ov M}_g(E/P,\,\beta _{[\sigma]})}$ 
containing $\sigma _0$ such that $X_0$ has the expected dimension. 
Now the proposition follows by the method of proof of the last part of the 
Theorem \ref{generic} and a variant of the Lemma \ref{largered} which 
states that if 
$\sigma _1$ is a $P$ reduction of $E$ and then there exists $P$ reduction 
$\sigma _2$ of $E$ such that $[\sigma _2]\le [\sigma _1]$ and that 
$[\sigma _2]$ satisfies the property ($*$) for $N_1$. This variant has a 
similar proof. This completes the proof of the Proposition \ref{main1}.
$\hfill \square$
\nocite{Harder1}
\nocite{SGA3}
\nocite{Lange}
\nocite{Bertram1}
\nocite{Bertram2}

\bibliographystyle{plain}

%\bibliography{master}

\end{document}